\documentclass[11pt,twoside]{article}
\usepackage{fancyhdr}
\usepackage[colorlinks,citecolor=blue,urlcolor=blue,linkcolor=blue,bookmarks=false]{hyperref}
\usepackage{amsfonts,epsfig,graphicx}
\usepackage{afterpage}
\usepackage{comment}
\usepackage{amsmath,amssymb,amsthm} 
\usepackage{fullpage}
\usepackage[T1]{fontenc} 
\usepackage{epsf} 
\usepackage{graphics} 
\usepackage{amsfonts,amsmath}
\usepackage[sort]{natbib} 
\usepackage{psfrag,xspace}
\usepackage{color,etoolbox}
\usepackage{subcaption} 
\usepackage{enumitem}

\usepackage{appendix}

\usepackage{tikz}
\usetikzlibrary{positioning,shapes.geometric}

\def\ind{\perp\!\!\!\perp}
\def\T{{ \mathrm{\scriptscriptstyle T} }}

\newcommand{\var}{\text{var}}

\newcommand{\Pb}{\mathbb{P}}
\newcommand{\Pn}{\mathbb{P}_n}
\newcommand{\Un}{\mathbb{U}_n}

\newcommand{\E}{\mathbb{E}}
\newcommand{\R}{\mathbb{R}}

\newcommand{\bias}{\text{bias}}
\newcommand{\Holder}{\text{H\"{o}lder}}

\DeclareMathOperator*{\argmin}{arg\,min}

\DeclareSymbolFont{bbold}{U}{bbold}{m}{n}
\DeclareSymbolFontAlphabet{\mathbbold}{bbold}
\newcommand{\one}{\mathbbold{1}}
\newtheorem{theorem}{Theorem}
\newtheorem{lemma}{Lemma}

\newtheorem{proposition}{Proposition}

\newtheorem{applemma}{Lemma}

\newtheorem{propositionapp}{Proposition}
\theoremstyle{definition}
\newtheorem{definition}{Definition}

\theoremstyle{remark}

\newtheorem{remark}{Remark}

\begin{document}

\def\spacingset#1{\renewcommand{\baselinestretch}%
{#1}\small\normalsize} \spacingset{1}

\raggedbottom
\allowdisplaybreaks[1]

%%%%%%%%%%%%%%%%%%%%%%%%%%%%%%%%%%%%%%%%%%%

  \title{\vspace*{-.4in} {Minimax rates for heterogeneous causal effect estimation}}
  \author{\\ Edward H.\ Kennedy$^1$, Sivaraman Balakrishnan$^{1,2}$, James M. Robins$^3$, Larry Wasserman$^{1,2}$ \\ \\
    $^1$Department of Statistics \& Data Science, Carnegie Mellon University \\ 
    $^2$Machine Learning Department, Carnegie Mellon University \\ 
    $^3$Departments of Biostatistics and Epidemiology, Harvard University \\ \\
    \texttt{\{edward, siva, larry\} @ stat.cmu.edu, robins@hsph.harvard.edu}  \\
\date{}
    }

  \maketitle
  \thispagestyle{empty}

\begin{abstract}
Estimation of heterogeneous causal effects -- i.e., how effects of policies and treatments vary across subjects -- is a fundamental task in causal inference. Many methods for estimating conditional average treatment effects (CATEs) have been proposed in recent years, but questions surrounding optimality have remained largely unanswered. In particular, a minimax theory of optimality has yet to be developed, with the minimax rate of convergence and construction of rate-optimal estimators remaining open problems.
In this paper we derive the minimax rate for CATE estimation, in a \Holder{}-smooth nonparametric model, and present a new local polynomial estimator, giving high-level conditions under which it is minimax optimal. 
Our minimax lower bound is derived via a localized version of the method of fuzzy hypotheses, combining lower bound constructions for nonparametric regression and functional estimation. Our proposed estimator can be viewed as a local polynomial R-Learner, based on a localized modification of higher-order influence function methods. 
The minimax rate we find exhibits several interesting features, including a non-standard elbow phenomenon and an unusual interpolation between nonparametric regression and functional estimation rates. The latter quantifies how the CATE, as an estimand, can be viewed as a regression/functional hybrid. 
\end{abstract}

\bigskip


\pagebreak

\section{Introduction}

In this paper we consider estimating the difference in regression functions
\begin{equation}
\tau(x) = \E(Y \mid X=x, A=1) - \E(Y \mid X=x, A=0)
\end{equation}
from an iid sample of observations of $Z=(X,A,Y)$. Let $Y^a$ denote the counterfactual outcome that would have been observed under treatment level $A=a$. 
Then, under the assumptions of consistency (i.e., $Y=Y^a$ if $A=a$), positivity (i.e., $\epsilon \leq \Pb(A = 1 \mid X) \leq 1-\epsilon$ with probability one, for some $\epsilon>0$), and no unmeasured confounding (i.e., $A \ind Y^a \mid X$), the quantity $\tau(x)$ also equals the conditional average treatment effect (CATE)
$$ \E(Y^1 - Y^0 \mid X=x) . $$ 
The CATE $\tau(x)$ gives a more individualized picture of treatment effects compared to the overall average treatment effect (ATE) $\E(Y^1-Y^0)$, and plays a crucial role in many fundamental tasks in causal inference, including assessing effect heterogeneity, constructing optimal treatment policies, generalizing treatment effects to new populations, finding subgroups with enhanced effects, and more. Further, these tasks have far-reaching implications across the sciences, from personalizing medicine to optimizing voter turnout. We refer to  \citet{hernan2020causal} (chapter 4), \citet{nie2017quasi}, \citet{kennedy2023towards}, and citations therein, for general discussion and review. \\

The simplest approach to CATE estimation would be to assume a low-dimensional parametric model for the outcome regression $\E(Y \mid X,A)$; then maximum likelihood estimates could be easily constructed, and under regularity conditions the resulting plug-in estimator would be minimax optimal. However, when $X$ has continuous components, it is typically difficult to specify a correct parametric model, and under misspecification the previously described approach could lead to substantial bias. This suggests the need for more flexible methods. Early work in flexible CATE estimation employed semiparametric models, for example partially linear models  assuming $\tau(x)$ to be constant, or structural nested models in which $\tau(x)$ followed some known parametric form, but leaving other parts of the distribution unspecified \citep{robinson1988root, robins1992estimating, robins1994correcting, van2003unified, van2006statistical, vansteelandt2014structural}. An important theme in this work is that the CATE can be much more structured and simple than the rest of the data-generating process. Specifically, the individual regression functions $\mu_a(x) = \E(Y \mid X=x,A=a)$ for each $a=0,1$ may be very complex (e.g., non-smooth or non-sparse), even when the difference $\tau(x) = \mu_1(x) - \mu_0(x)$ is very smooth or sparse, or even constant or zero. We refer to \citet{kennedy2023towards} for some recent discussion of this point. \\

%To the best of our knowledge, \citet{robins2008higher} (Section 5.2) seems to be the first proposed nonparametric model where the CATE has its own smoothness, possibly very different from that of the individual regression functions; however this model considered experiments with known propensity scores.

More recently there has been increased emphasis on incorporating nonparametrics and machine learning tools for CATE estimation. We  briefly detail two especially relevant streams of this recent literature, based on so-called DR-Learner and R-Learner methods, both of which rely on doubly robust-style estimation. The DR-Learner is a model-free meta-algorithm first proposed by \citet{van2006statistical} (Section 4.2), which essentially takes the components of the classic doubly robust estimator of the ATE, and rather than averaging, instead regresses on covariates. It has since been specialized to particular methods, e.g., cross-validated ensembles \citep{luedtke2016super}, kernel \citep{lee2017doubly, fan2019estimation, zimmert2019nonparametric} and series methods \citep{semenova2017estimation}, empirical risk minimization \citep{foster2019orthogonal}, and linear smoothers \citep{kennedy2023towards}.  On the other hand, the R-Learner is a flexible adaptation of the double-residual regression method originally built for partially linear models \citep{robinson1988root}, with the first nonparametric version proposed by \citet{robins2008higher} (Section 5.2) using series methods. The R-Learner has since  been adapted to RKHS regression \citep{nie2017quasi}, lasso \citep{zhao2017selective, chernozhukov2017orthogonal}, and local polynomials \citep{kennedy2023towards}.  Many flexible non-doubly robust methods have also been proposed in recent years, often based on inverse-weighting or direct regression estimation \citep{foster2011subgroup, imai2013estimating, athey2016recursive, shalit2017estimating, wager2018estimation, kunzel2019metalearners, hahn2020bayesian}. \\

Despite the wide variety of methods available for flexible CATE estimation, questions of optimality have remained mostly unsolved. \citet{gao2020minimax} studied minimax optimality, but in a specialized model where the propensity score has zero smoothness, and covariates are non-random; this model does not reflect the kinds of assumptions typically used in practice, e.g., in the papers cited in the previous paragraph.  Some but not all of these papers derive upper bounds on the error of their proposed CATE estimators; in the best case, these take the form of an oracle error rate (which would remain even if the potential outcomes $(Y^1-Y^0)$ were observed and regressed on covariates), plus some contribution coming from having to estimate nuisance functions (i.e., outcome regressions and propensity scores). The fastest rates we are aware of come from \citet{foster2019orthogonal} and \citet{kennedy2023towards}. \citet{foster2019orthogonal} studied global error rates, obtaining an oracle error plus sums of squared $L_4$ errors in all nuisance components. \citet{kennedy2023towards} studied pointwise error rates, giving two main results; in the first, they obtain the oracle error plus a product of nuisance errors, while in the second, they obtain a faster rate via undersmoothing (described in more detail in Section \ref{sec:finalrate}). However, since these are all upper bounds on the errors of particular procedures, it is unknown whether these rates are optimal in any sense, and if they are not, how they might be improved upon. In this paper we resolve these questions (via the minimax framework, in a nonparametric model that allows components of the data-generating process to be infinite-dimensional, yet smooth in the \Holder{} sense). \\

More specifically, in Section \ref{sec:lowerbound}  we derive a lower bound on the minimax rate of CATE estimation, indicating the best possible (worst-case) performance of any estimator, in a model where the CATE, regression function, and propensity score are \Holder{}-smooth functions.  Our derivation uses an adaptation of the method of fuzzy hypotheses, which is specially localized compared to the constructions previously used for obtaining lower bounds in functional estimation and hypothesis testing \citep{ibragimov1987some, birge1995estimation, nemirovski2000topics, ingster2003nonparametric, tsybakov2009introduction, robins2009semiparametric}. In Section \ref{sec:upperbound}, we confirm that our minimax lower bound is tight (under some conditions), by proposing and analyzing a new local polynomial R-Learner, using localized adaptations of higher order influence function methodology \citep{robins2008higher, robins2009quadratic, robins2017minimax}. In addition to giving a new estimator that is provably optimal (under some conditions, e.g., on how well the covariate density is estimated), our results also confirm that previously proposed estimators were not generally optimal in this smooth nonparametric model. Our minimax rate also sheds light on the nature of the CATE as a statistical quantity, showing how it acts as a regression/functional hybrid: for example, the rate interpolates between nonparametric regression and functional estimation, depending on the relative smoothness of the CATE and nuisance functions (outcome regression and propensity score). \\

\section{Setup \& Notation}

We consider an iid sample of $n$ observations of $Z=(X,A,Y)$ from distribution $\Pb$, where $X \in [0,1]^d$ denotes covariates, $A \in \{0,1\}$  a treatment or policy indicator, and $Y \in \R$ an outcome of interest. We let $F$ denote the distribution function of the covariate $X$ (with density $f$ as needed), and let 
\begin{align*}
\pi(x) &= \Pb(A=1 \mid X=x) \\
\eta(x) &= \E(Y \mid X=x) \\
\mu_a(x) &= \E(Y \mid X=x, A=a)
\end{align*}
denote the propensity score, marginal, and treatment-specific outcome regressions, respectively. We sometimes omit arguments from functions to ease notation, e.g., note that $\tau=(\eta-\mu_0)/\pi$. We also index quantities by a distribution $P$ when needed, e.g., $\tau(x)$ under a particular distribution $P$ is written $\tau_P(x)$; depending on context, no indexing means the quantity is evaluated at the true $\Pb$, e.g., $\tau(x)=\tau_\Pb(x)$.  \\

Our goal is to study estimation of the CATE $\tau(x)=\mu_1(x)-\mu_0(x)$ at a point $x_0 \in (0,1)^d$, with error quantified by mean absolute error 
$$ \E \left| \widehat\tau(x_0) - \tau(x_0) \right| . $$
As detailed in subsequent sections, we work  in a nonparametric model $\mathcal{P}$ whose components are infinite-dimensional functions but with some smoothness. 
We say a function is $s$-smooth if it belongs to a \Holder{} class with index $s$; this  essentially means it has $s-1$ bounded derivatives, and the highest order derivative is continuous. To be more precise, let $\lfloor s \rfloor$ denote the largest integer strictly smaller than $s$, and let $D^\alpha = \frac{\partial^\alpha }{\partial x_1^{\alpha_1} \dots \partial x_d^{\alpha_d}}$ denote the partial derivative operator. Then the \Holder{} class with index $s$ contains all functions $g: \mathcal{X} \rightarrow \R$ that are $\lfloor s \rfloor$ times continuously differentiable, with derivatives up to order $\lfloor s \rfloor$ bounded, i.e., 
$$  \left| D^\alpha g(x) \right| \leq C < \infty $$
for all $\alpha=(\alpha_1,\dots,\alpha_d)$ with $\sum_j \alpha_j \leq \lfloor s \rfloor$ and for all $x \in \mathcal{X}$, and with $\lfloor s \rfloor$-order derivatives satisyfing the Lipschitz condition
$$ \left| D^\beta g(x) - D^\beta g(x') \right| \leq C \| x - x' \|^{s-\lfloor s \rfloor}  $$
for some $C<\infty$, for all  $\beta=(\beta_1,\dots,\beta_d)$ with $\sum_j \beta_j = \lfloor s \rfloor$, and for all $x,x' \in \mathcal{X}$, where for a vector $v \in \R^d$ we let $\|v \|$ denote the Euclidean norm.  Sometimes \Holder{} classes are referenced by both the smoothness $s$ and constant $C$, but we focus our discussion on the smoothness $s$ and omit the constant, which is assumed finite and independent of $n$. \\

We write the squared $L_2(Q)$ norm of a function as $\| g \|_{Q}^2 = { \int g(z)^2 \ dQ(z)}$. The sup-norm is denoted by $\|f \|_\infty = \sup_{z \in \mathcal{Z}} | f(z) |$. For a matrix $A$ we let $\| A \|$ and $\| A \|_2$ denote the operator/spectral and Frobenius norms, respectively, and let $\lambda_{\min}(A)$ and $\lambda_{\max}(A)$ denote the minimum and maximum eigenvalues of $A$, respectively. We write $a_n \lesssim b_n$ if $a_n \leq C b_n$ for $C$ a positive constant independent of $n$, and $a_n \asymp b_n$ if $a_n \leq C b_n$ and $b_n \leq C a_n$ (i.e., if $a_n \lesssim b_n$ and $b_n \lesssim a_n$).  We write $a_n \sim b_n$ to mean that $a_n$ and $b_n$ are proportional, i.e., $a_n= C b_n$ for some $C$. We also use $a \vee b = \max(a,b)$ and $a \wedge b = \min(a,b)$. We use the shorthand $\Pn(f) = \Pn\{f(Z)\} = \frac{1}{n} \sum_{i=1}^n f(Z_i)$ to write sample averages, and similarly $\Un(f) = \Un\{f(Z_1,Z_2)\}=\frac{1}{n(n-1)} \sum_{i \neq j} f(Z_i,Z_j)$ for the U-statistic measure. \\

\section{Fundamental Limits}
\label{sec:lowerbound}

In this section we derive a lower bound on the minimax rate for CATE estimation. This result has several crucial implications, both practical and theoretical. First, it gives a benchmark for the best possible performance of any CATE estimator in the nonparametric model defined in Theorem \ref{thm:lowerbound}. In particular, if an estimator is shown to attain this benchmark, then one can safely conclude the estimator cannot be improved, at least in terms of worst-case rates, without adding assumptions; conversely, if the benchmark is \emph{not} shown to be attained, then one should continue searching for other better estimators (or better lower or upper risk bounds).  Second,  a tight minimax lower bound is important in its own right as a measure of the fundamental limits of CATE estimation, illustrating precisely how difficult CATE estimation is in a statistical sense.  The main result of this section is given in Theorem \ref{thm:lowerbound} below. It is finally proved and discussed in detail in Section \ref{sec:finalrate}. \\

\begin{theorem} \label{thm:lowerbound}
For $x_0 \in (0,1)^d$, let $\mathcal{P}$ denote the model where:
\begin{enumerate}
\item $f(x)$ is bounded above by a constant, 
\item $\pi(x)$ is $\alpha$-smooth, %and $\epsilon \leq \pi(x) \leq 1-\epsilon$ for some $\epsilon>0$, 
\item $\mu_0(x)$ is $\beta$-smooth, and
\item $\tau(x)$ is $\gamma$-smooth.
\end{enumerate}
Let $s \equiv (\alpha+\beta)/2$. 
Then for $n$ larger than a constant depending on $(\alpha,\beta,\gamma,d)$, the minimax rate is lower bounded as
\begin{align*}
\inf_{\widehat\tau} \sup_{P \in \mathcal{P}} \E_P | \widehat\tau(x_0) - \tau_P(x_0)| \gtrsim \begin{cases}
n^{-1/\left(1 + \frac{d}{2\gamma} + \frac{d}{4s} \right) } & \text{ if } s < \frac{d/4}{1+d/2\gamma}  \\
n^{-1/\left(2 + \frac{d}{\gamma}  \right) } & \text{ otherwise}  . 
\end{cases}
\end{align*}
\end{theorem}

\bigskip

\begin{remark}
In Appendix \ref{sec:model2} we also give results (both lower and upper bounds) for the model that puts smoothness assumptions on the marginal regression $\eta(x)=\E(Y \mid X=x)$, instead of the control regression $\mu_0(x)=\E(Y \mid X=x,A=0)$. Interestingly, the minimax rates differ in these two models, but only in the regime  where the regression function is more smooth than the propensity score (i.e., $\beta > \alpha$). Specifically, when $\eta$ is $\beta$-smooth, the minimax rate from Theorem \ref{thm:lowerbound} holds but with $s=(\alpha+\beta)/2$ replaced by $\min(\alpha,s)$.  \\
\end{remark}

Crucially, Condition 4 allows the CATE $\tau(x)$ to have its own smoothness  $\gamma$, which is necessarily at least the regression smoothness $\beta$, but can also be much larger, as described in the Introduction.  
We defer discussion of the details of the overall minimax rate of Theorem \ref{thm:lowerbound} to Section \ref{sec:finalrate}, moving first to a proof of the result.  \\

\begin{remark}
For simplicity,  the lower bound result in Theorem \ref{thm:lowerbound} is given for a large model in which the covariate density is only bounded. However,  as discussed in detail later in Section \ref{sec:upperbound} and Remark \ref{rem:cond1}, for the stated rate to be attainable more conditions on the covariate density are required. In Section \ref{sec:appcovdensity} of the Appendix, we give a particular submodel of $\mathcal{P}$ under which upper and lower bounds on the minimax rate match up to constants; it will be important in future work to further elucidate the role of the covariate density in CATE estimation.\\
\end{remark}

The primary strategy in deriving minimax lower bounds is to construct distributions that are similar enough that they are statistically indistinguishable, but for which the parameter of interest is maximally separated; this implies no  estimator can have error uniformly smaller than this separation. More specifically, we derive our lower bound using a localized version of the method of fuzzy hypotheses \citep{ibragimov1987some, birge1995estimation, nemirovski2000topics, ingster2003nonparametric, tsybakov2009introduction, robins2009semiparametric}. In the classic Le Cam two-point method, which can be used to derive minimax lower bounds for nonparametric regression at a point \citep{tsybakov2009introduction}, it suffices to consider a pair of distributions that differ locally; however, for nonlinear functional estimation, such pairs give bounds that are too loose. One instead needs to construct pairs of \emph{mixture} distributions, which can be viewed via a prior over distributions in the model \citep{birge1995estimation, tsybakov2009introduction, robins2009semiparametric}. Our construction combines these two approaches via a localized mixture, as will be described in detail in the next subsection. \\

\begin{remark} \label{rem:lowsmooth}
In what follows we focus on the lower bound in the low smoothness regime where $s < \frac{d/4}{1+d/2\gamma}$. The $n^{-1/(2+d/\gamma)}$ lower bound for the high smoothness regime matches the classic smooth nonparametric regression rate, and follows from a standard two-point argument, using the same construction as in Section 2.5 of  \citet{tsybakov2009introduction}. \\
\end{remark}

The following lemma, adapted from Section 2.7.4 of \citet{tsybakov2009introduction}, provides the foundation for the minimax lower bound result of this section. \\

\begin{lemma}[\citet{tsybakov2009introduction}] \label{lem:minimax}
Let $P_\lambda$ and $Q_\lambda$ denote distributions in $\mathcal{P}$ indexed by a vector $\lambda=(\lambda_1,\dots,\lambda_k)$, with $n$-fold products denoted by $P_\lambda^n$ and $Q_\lambda^n$, respectively.  
Let $\varpi$ denote a prior distribution over $\lambda$.  
If 
$$ H^2\left(\int P_\lambda^n \ d\varpi(\lambda), \int Q_\lambda^n \ d\varpi(\lambda) \right) \leq \alpha < 2$$ 
and 
$$ \psi(P_\lambda) - \psi(Q_{\lambda'}) \geq  s >0 $$
for a functional $\psi: \mathcal{P} \mapsto \R$ and for all $\lambda,\lambda'$, 
then
$$ \inf_{\widehat\psi} \sup_{P \in \mathcal{P}} \E_P \left\{ \ell\left( \left| \widehat\psi - \psi(P)\right| \right) \right\} \geq \ell(s/2) \left( \frac{1 - \sqrt{\alpha(1-\alpha/4)}}{2} \right) $$
for any monotonic non-negative loss function $\ell$.
\end{lemma}

\bigskip

Lemma \ref{lem:minimax} illuminates the three ingredients for deriving a minimax lower bound, and shows how they interact. The ingredients are: (i) a pair of mixture distributions, (ii) the distance between their $n$-fold products, which is ideally small, and  (iii) the separation of the parameter of interest under the mixtures, which is ideally large. Finding the right minimax lower bound requires balancing these three ingredients appropriately: with too much distance or not enough separation,  the lower bound will be too loose. In the following subsections we describe these three ingredients in detail. \\

\begin{comment}
\begin{proof}
For any monotonic non-negative loss function $\ell$ we have
\begin{align*}
 \inf_{\widehat\psi} \sup_{P \in \mathcal{P}} \E_P \left\{ \ell\left( | \widehat\psi - \psi(P)| \right) \right\} &\geq \ell(s/2) \inf_{\widehat\psi} \sup_{P \in \mathcal{P}} P\left( | \widehat\psi - \psi(P)| \geq s/2 \right)  \\
&\geq \ell(s/2) \left( \frac{ 1 - \sqrt{\alpha(1-\alpha/4)} }{2} \right)
 \end{align*}
 where the first line follows by Markov's inequality, and the second from Theorem 2.15(ii) of \citet{tsybakov2009introduction} after taking  $\beta_0=\beta_1=0$.
\end{proof}
\end{comment}

\subsection{Construction}

In this subsection we detail the distributions $P_\lambda$ and $Q_\lambda$ used to construct the minimax lower bound. 
The main idea is to mix constructions for nonparametric regression and functional estimation, by perturbing the CATE with a bump at the point $x_0$, and to also use a mixture of perturbations of the propensity score and regression functions $\pi$ and $\mu_0$,  locally near $x_0$.  \\

For our lower bound results, we work in the setting where $Y$ is binary; this is mostly to ease notation and calculations. Note however that this still yields a valid lower bound in the general continuous $Y$ case, since a lower bound in the strict submodel where $Y$ is binary is also a lower bound across the larger model $\mathcal{P}$. Importantly, when $Y$ is binary, the density $p$ of an observation $Z$ can be  indexed via either the quadruple $(f,\pi,\mu_0,\mu_1)$, or $(f,\pi,\mu_0,\tau)$; here we make use of the latter parametrization (and in the appendix we consider the $(f,\pi,\eta,\tau)$ parametrization). 
We first give the construction for the $\alpha \geq \beta$ case in the definition below, and then go on to discuss details (and in Appendix \ref{sec:model2} we give constructions for all other regimes). \\

\begin{definition}[Distributions $P_\lambda$ and $Q_\lambda$]  \label{def:construction}
Let:
\begin{enumerate}
\item $B: \R^d \rightarrow \R$ denote a $C^\infty$ function with $B(x)=1$ for $x \in [-1/4, 1/4]^d$, and $B(x)=0$ for $x \notin [-1/2,1/2]^d$,
\item $\mathcal{C}_h(x_0)$ denote the cube centered at $x_0 \in (0,1)^d$ with sides of length $h \leq 1/4$, 
\item $(\mathcal{X}_1,\dots,\mathcal{X}_k)$ denote a partition of $\mathcal{C}_h(x_0)$ into $k$ cubes of equal size (for $k$ an integer raised to the power $d$), with midpoints $(m_1,\dots,m_k)$, so each cube $\mathcal{X}_j = \mathcal{C}_{h/k^{1/d}}(m_j)$ has side length $h/k^{1/d}$.
\end{enumerate}
Then for $\lambda_j \in \{-1,1\}$ define the functions
\begin{align*}
\tau_h(x) &= h^\gamma B\left( \frac{x-x_0}{2h} \right)  \\
\mu_{0\lambda}(x) &= \frac{1}{2} +  \left( h/k^{1/d} \right)^{\beta} \sum_{j=1}^k \lambda_j B \left( \frac{x-m_j}{h/k^{1/d}}  \right) \\
\pi_\lambda(x) &=   \frac{1}{2} +   \left(  h/k^{1/d} \right)^{\alpha} \sum_{j=1}^k \lambda_j B \left( \frac{x-m_j}{h/k^{1/d}}  \right) \\
 f(x) &= \one(x \in \mathcal{S}_{hk}) \Big/ \left\{ 1  - \left( \frac{4^d - 1}{2^d} \right) h^d \right\}
\end{align*}
where $\mathcal{S}_{hk} = \left\{ \bigcup_{j=1}^k \mathcal{C}_{h/2k^{1/d}}(m_j) \right\} \bigcup \left\{ [0,1]^d \setminus \mathcal{C}_{2h}(x_0) \right\}$. 
Finally let the distributions $P_\lambda$ and $Q_\lambda$ be defined via the densities
\begin{align*}
p_\lambda &= (f,1/2, \mu_{0\lambda}-\tau_h/2,\tau_h) \\
q_\lambda &= (f,\pi_\lambda,\mu_{0\lambda},0) . 
\end{align*}
\end{definition}

\bigskip

\begin{remark}
Under the $(f,\pi,\eta,\tau)$ parametrization, since $\eta = \pi\tau+\mu_0$, the above densities can equivalently be written as $p_\lambda=(f,1/2, \mu_{0\lambda},\tau_h)$ and $q_\lambda=(f,\pi_\lambda,\mu_{0\lambda},0)$. \\
\end{remark}

Figure \ref{fig:construction} shows an illustration of our construction in the $d=1$ case. As mentioned above,  the CATE is perturbed with a bump at $x_0$ and the nuisance functions $\pi$ and $\mu_0$ with bumps locally near $x_0$. The regression function $\mu_0$ is perturbed under both $P_\lambda$ and $Q_\lambda$, since it is less smooth than the propensity score in the $\alpha \geq \beta$ case. The choices of the CATE mimic those in the two-point proof of the lower bound for nonparametric regression at a point (see, e.g., Section 2.5 of \citet{tsybakov2009introduction}), albeit with a particular flat-top bump function, while the choices of nuisance functions $\pi$ and $\mu_0$ are more similar to those in the lower bound for functionals such as the expected conditional covariance (cf.\ Section 4 of \citet{robins2009semiparametric}). In this sense our construction can be viewed as combining those for nonparametric regression and functional estimation, similar to  \citet{shen2020optimal}. In what follows we remark on some important details.  \\

\begin{remark}
Section 3.2 of \citet{shen2020optimal} used a similar construction for conditional variance estimation. Some important distinctions are: (i) they focused on the univariate and low smoothness setting; (ii) in that problem there is only one nuisance function, so the null can be a point rather than a mixture distribution; and (iii) they use a different, arguably more complicated, approach to bound the distance between distributions. Our work can thus be used to generalize such variance estimation results to arbitrary dimension and smoothness. \\
\end{remark}

\begin{figure}[h!] \label{fig:construction}
\begin{subfigure}{.49\textwidth}
  \centering
  % include first image
  \includegraphics[width=1.02\linewidth]{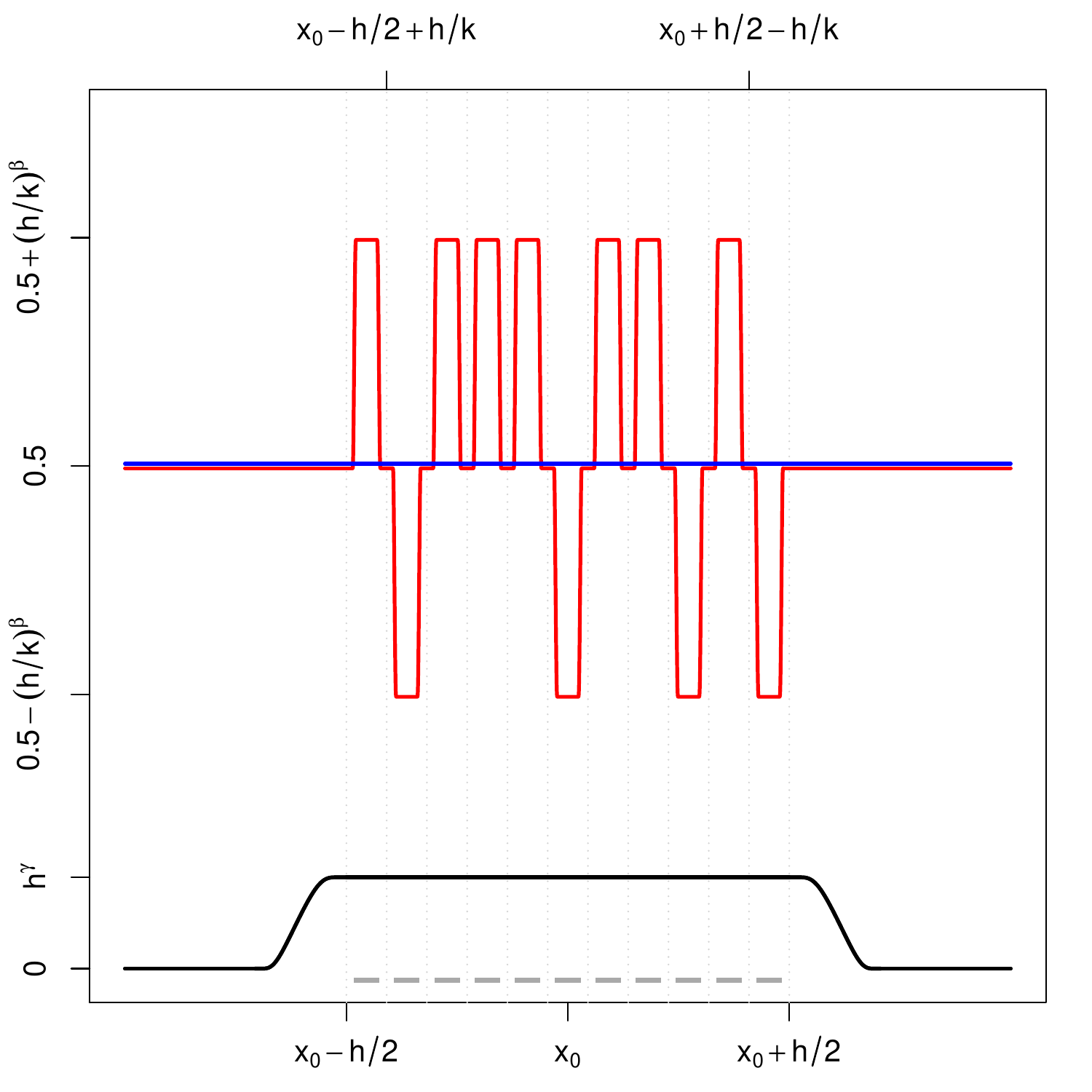}  
  \caption{Null $P_\lambda$}
  \label{fig:null}
\end{subfigure}
\begin{subfigure}{.49\textwidth}
  \centering
  % include second image
  \includegraphics[width=1.02\linewidth]{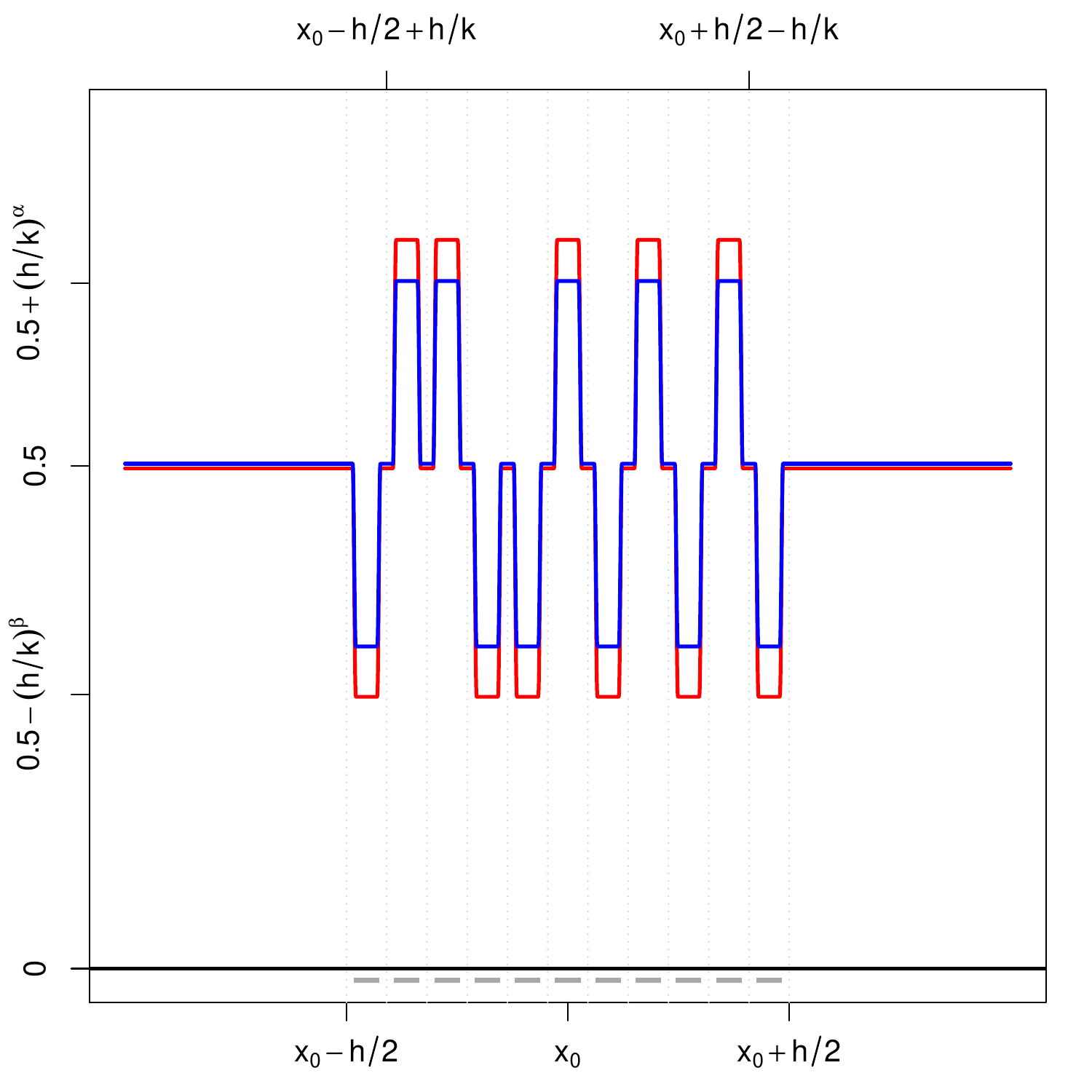}  
  \caption{Alternative $Q_\lambda$}
  \label{fig:alt}
\end{subfigure}
\caption{Minimax lower bound construction in $d=1$ case. An example null density $p_\lambda$ is displayed in panel (a) and an alternative density $q_\lambda$ in panel (b). The black, red, and blue lines denote the CATE, \emph{marginal} outcome regression $\eta=\pi\tau+\mu_0$, and propensity score functions, respectively, and the gray line denotes the support of the covariate density.  }
\label{fig:construction}
\end{figure}

\bigskip

First we remark on the choice of CATE in the construction. As mentioned above, the bump construction resembles that of the standard Le Cam lower bound for nonparametric regression at a point, but differs in that we use a specialized  bump function with a flat top. Crucially, this choice ensures the CATE  is constant and equal to $h^\gamma$ for all $x$ in the cube $\mathcal{C}_h(x_0)$ centered at $x_0$ with sides of length $h$, and  that it is equal to zero for all $x \notin \mathcal{C}_{2h}(x_0)$, i.e., outside the cube  centered at $x_0$ with side length $2h$. The fact that the CATE is constant  across the top of the bump (which will be the only place where observations appear near $x_0$) eases Hellinger distance calculations substantially. It is straightforward to check that the CATE $\tau_h(x)$ is $\gamma$-smooth in this construction (see page 93 of \citet{tsybakov2009introduction}).   \\

\begin{remark}
One example of a bump function $B$ satisfying the conditions above is 
$$ B(x) =  \left[ 1 + \frac{ g\{(4x)^2-1\} }{ g\{2-(4x)^2\} } \right]^{-1} $$
for $g(t) = \exp(-1/t) \one(t>0)$. \\
\end{remark}

For the propensity score and regression functions, we similarly have
\begin{align*}
B \left( \frac{x-m_j}{h/k^{1/d}}  \right) = \begin{cases}
1 & \text{ for } x \in \mathcal{C}_{h/2k^{1/d}}(m_j) \\
0 & \text{ for } x \notin \mathcal{C}_{h/k^{1/d}}(m_j) 
\end{cases}
\end{align*}
i.e., each bump equals one on the half-$h/k^{1/d}$ cube around $m_j$, and is identically zero outside the main  larger $h/k^{1/d}$ cube around $m_j$. It is again straightforward to check that $\pi_\lambda(x)$ and $\mu_{0\lambda}(x)$ are $\alpha$- and $\beta$-smooth, respectively. \\
 
The covariate density is chosen to be uniform, but on the set $\mathcal{S}_{hk}$  that captures the middle of all the nuisance bumps $\left\{ \bigcup_{j=1}^k \mathcal{C}_{h/2k^{1/d}}(m_j) \right\}$, together with the space $\left\{ [0,1]^d \setminus \mathcal{C}_{2h}(x_0) \right\}$ away from $x_0$. Importantly, this choice ensures there is only mass where the nuisance bumps $B \left( \frac{x-m_j}{h/k^{1/d}}  \right)$ are constant and non-zero (and where $\tau_h(x)=h^\gamma$), or else far away from $x_0$, where the densities are the same under $P_\lambda$ and $Q_\lambda$. Note that, as $h \rightarrow 0$, the Lebesgue measure of the set $\mathcal{S}_{hk}$ tends to one, and the covariate density tends towards a standard uniform distribution.  \\

\subsection{Hellinger Distance}

As mentioned previously, deriving a tight minimax lower bound requires carefully balancing the distance between distributions in our construction. To this end, in this subsection we bound the Hellinger distance between the $n$-fold product mixtures $\int P_\lambda^n \ d\varpi(\lambda)$ and $\int Q_\lambda^n \ d\varpi(\lambda)$, for $\varpi$ a uniform prior distribution, so that $(\lambda_1,\dots,\lambda_k)$ are iid Rademachers. \\

In general these product densities can be complicated, making direct distance calculations difficult. However, Theorem 2.1 from \citet{robins2009semiparametric} can be used to relate the distance between the $n$-fold products to those of simpler posteriors over a single observation. In the following lemma we adapt this result to localized constructions like those in Definition \ref{def:construction}. \\

\begin{lemma} \label{lem:hellbound}
Let $P_\lambda$ and $Q_\lambda$ denote distributions indexed by a vector $\lambda=(\lambda_1,\dots,\lambda_k)$, and let $\mathcal{Z}=\cup_{j=1}^k \mathcal{Z}_j$ denote a partition of the sample space. 
Assume:
\begin{enumerate}
\item $P_\lambda(\mathcal{Z}_j)=Q_\lambda(\mathcal{Z}_j)=p_j$ for all $\lambda$, and
\item the conditional distributions $\one_{\mathcal{Z}_j} dP_\lambda / p_j$ and $\one_{\mathcal{Z}_j} dQ_\lambda / p_j$ (given an observation is in $\mathcal{Z}_j$) do not depend on $\lambda_\ell$ for $\ell \neq j$, and only differ on partitions $j \in S \subseteq \{1,\dots,k\}$.
\end{enumerate}
For a prior distribution $\varpi$ over $\lambda$, let $\overline{p}=\int p_\lambda \ d\varpi(\lambda)$ and $\overline{q}=\int q_\lambda \ d\varpi(\lambda)$, and define 
\begin{align*}
\delta_1 &= \max_{j \in S} \sup_\lambda \int_{\mathcal{Z}_j} \frac{(p_\lambda - \overline{p})^2}{p_\lambda p_j} \ d\nu \nonumber \\
\delta_2 &= \max_{j \in S} \sup_\lambda \int_{\mathcal{Z}_j} \frac{(q_\lambda - p_\lambda)^2}{p_\lambda p_j} \ d\nu \label{eq:bdcomps} \\
\delta_3 &= \max_{j \in S} \sup_\lambda \int_{\mathcal{Z}_j} \frac{(\overline{q} - \overline{p})^2}{p_\lambda p_j} \ d\nu \nonumber
\end{align*}
for  a dominating measure $\nu$. 
If $\overline{p}/p_\lambda \leq b < \infty$ and $n p_j \max(1,\delta_1,\delta_2) \leq b$ for all $j$, 
then $H^2\left(\int P_\lambda^n \ d\varpi(\lambda), \int Q_\lambda^n \ d\varpi(\lambda) \right) $ is bounded above by
$$  Cn \left( \sum_{j \in S} p_j \right) \left\{ n \left( \max_{j \in S} p_j \right) \Big(\delta_1 \delta_2 + \delta_2^2 \Big) + \delta_3 \right\} $$ 
for a constant $C$ only depending on $b$.
\end{lemma}

\bigskip

In the next proposition, we bound the quantities from Lemma  \ref{lem:hellbound} and put the results together to obtain a bound on the desired Hellinger distance between product mixtures. \\

\begin{proposition} \label{prop:hellbound}
Assume $h \leq 1/4$ and $h^\gamma + 2(h/k^{1/d})^\beta \leq 1 - 4 \varepsilon$ for some $\varepsilon \in (0,1/4)$, and take $h^\gamma = 4(h/k^{1/d})^{2s}$ for $s \equiv (\alpha+\beta)/2$. Then for the distributions $P_\lambda$ and $Q_\lambda$ from Definition \ref{def:construction}, with  $\varpi$ the uniform distribution over $\{-1,1\}^k$,  the quantities $\delta_j$ from Lemma \ref{lem:hellbound} satisfy
\begin{align*}
\delta_1 \leq \left( \frac{2^{d+1} \|B \|_2^2}{\varepsilon} \right) \left( h/k^{1/d} \right)^{2\beta} , \ \ \  \delta_2 \leq \left( \frac{2^{d+1} \|B \|_2^2}{\varepsilon} \right) \left( h/k^{1/d} \right)^{2\alpha} , \ \ \ \delta_3 = 0 , 
\end{align*}
and $ p_j \leq 2 (h/2)^d/k$. Further $H^2\left(\int P_\lambda^n \ d\varpi(\lambda), \int Q_\lambda^n \ d\varpi(\lambda) \right) $ is bounded above by
$$ 4C \left( \frac{ 2 \|B \|_2^2}{\varepsilon} \right)^2 \left( \frac{n^2 h^{2d}}{k} \right) \left\{ \left( h/k^{1/d} \right)^{4s} + \left( h/k^{1/d} \right)^{4\alpha} \right\} $$ 
for $C$ a constant only depending on $\varepsilon$.
\end{proposition}

\bigskip

Before moving to the proof of Proposition \ref{prop:hellbound}, we briefly discuss and give some remarks. Compared to the Hellinger distance arising in the average treatment effect or expected conditional covariance lower bounds \citep{robins2009semiparametric}, there is an extra $h^d$ factor in the numerator. Of course, one cannot simply repeat those calculations with $k/h^d$ bins, since then for example the $k^{-4s/d}$ term would also be inflated to $(k/h^d)^{-4s/d}$;  our carefully localized construction is crucial to obtain the right rate in this case. We also note that the choice $h^\gamma = 4(h/k^{1/d})^{2s}$ is required for ensuring that the averaged densities $\overline{p}(z)$ and $\overline{q}(z)$ are equal (implying that $\delta_3=0$); specifically this equalizes the CATE bump under $P_\lambda$ with the squared nuisance bumps under $Q_\lambda$. \\

%We also note that the conditions $h^\gamma \leq 1/4$ and $h^\gamma + 2k^{-\beta/d} \leq 1 - 4 \epsilon$

\begin{proof}

Here the relevant partition of the sample space $\mathcal{X} \times \mathcal{A} \times \mathcal{Y} = [0,1]^d \times \{0,1\} \times \{0,1\}$ is  $\mathcal{Z}_j = \mathcal{C}_{h/2k^{1/d}}(m_j) \times \{0,1\} \times \{0,1\}$, $j=1,\dots,k$, along with $\mathcal{Z}_j'$, which partitions the space $[0,1]^d / \mathcal{C}_{2h}(x_0)$ away from $x_0$ into disjoint cubes with side lengths $h/2k^{1/d}$. 
Thus we have
$$ P_\lambda(\mathcal{Z}_j) = P_\lambda(\mathcal{Z}_j') = Q_\lambda(\mathcal{Z}_j) = Q_\lambda(\mathcal{Z}_j') = p_j $$
where $p_j =  \int \one\{x \in C_{h/2k^{1/d}}(m_j)\} f(x) \ dx$. In Appendix Section \ref{app:prophellbound} we show that $(h/2)^d/k \leq p_j \leq 2 (h/2)^d/k$ when $h \leq 1/4$, and so is proportional to the volume of a cube with side lengths $h/2k^{1/d}$. 
Further   the conditional distributions $\one_{\mathcal{Z}_j} dP_\lambda / p_j$ and $\one_{\mathcal{Z}_j} dQ_\lambda / p_j$ do not depend on $\lambda_\ell$ for $\ell \neq j$, since $\lambda_j$ only changes the density in $\mathcal{Z}_j$.  
In what follows we focus on the partitions $\mathcal{Z}_j$ and not $\mathcal{Z}_j'$, since the distributions are exactly equal on the latter (in the language of Lemma \ref{prop:hellbound}, the set $S$ indexes only the partitions $\mathcal{Z}_j$, and note that for this set we have $\sum_{j \in S} p_j=kp_j \leq 2 (h/2)^d$).  
 Note when $(\lambda_1,\dots,\lambda_k)$ are iid Rademacher random variables the marginalized densities from Lemma \ref{lem:hellbound} are
\begin{align*}
\overline{p}(z) &= f(x) \left\{ \frac{1}{4} +   (2a-1)(2y-1)  \frac{h^\gamma}{4} B\left( \frac{x-x_0}{2h} \right) \right\} \\
\overline{q}(z) &= f(x) \left\{ \frac{1}{4}  +  (2a-1)(2y-1)  \left( h/k^{1/d} \right)^{2s} \sum_{j=1}^k B \left( \frac{x-m_j}{h/k^{1/d}}  \right) ^2 \right\} .
\end{align*}

\bigskip

The first step is to show that relevant densities and density ratios are appropriately bounded. We give these details in Appendix Section \ref{app:prophellbound}. Next it remains to bound the quantities $\delta_1$, $\delta_2$, and $\delta_3$. \\

%$$ h^\gamma + 2k^{-\beta/d} \leq 1-4\epsilon  \implies  \epsilon \leq 1/4 - \left( \frac{ k^{-\beta/d} }{2}  + \frac{h^\gamma}{4} \right)  $$
%$$ h \leq 1/4 \implies  1  - \left( \frac{4^d - 1}{2^d} \right) h^d  \geq 1  - \left( \frac{4^d - 1}{4^d 2^d} \right) \geq 1 - \frac{1}{2^d} \geq \frac{1}{2} $$
%$$  \frac{\overline{p}(z)}{p_\lambda(z)} \leq \frac{\frac{1}{4} + \frac{h^\gamma}{4}}{\frac{1}{4} - \frac{ k^{-\beta/d} }{2}  - \frac{h^\gamma}{4}  } \leq \epsilon^{-1} \left( \frac{1}{4} + \frac{1}{4^{\gamma+1}} \right) \leq \frac{1}{2\epsilon} $$

We begin with $\delta_3$, the distance between marginalized densities $\overline{p}$ and $\overline{q}$, which is tackled somewhat differently from $\delta_1$ and $\delta_2$. Because we take $(h/k^{1/d})^{2s} = h^\gamma/4$ it follows that $\overline{q}(z) - \overline{p}(z)$ equals
\begin{equation}
 (2a-1)(2y-1) f(x) \left\{   \left( h/k^{1/d} \right)^{2s} \sum_{j=1}^k B \left( \frac{x-m_j}{h/k^{1/d}}  \right)^2 - \frac{h^\gamma}{4} B\left( \frac{x-x_0}{2h} \right) \right\} = 0  , \label{eq:delta3}
\end{equation}
since  $f(x) = 0$ for $x \notin \mathcal{S}_{hk}$ and
\begin{align*}
B \left( \frac{x-m_j}{h/k^{1/d}}  \right)=B\left(\frac{x-x_0}{2h} \right) = 0 & \ \text{ for } x \in \left\{ [0,1]^d \setminus \mathcal{C}_{2h}(x_0) \right\} \\
B \left( \frac{x-m_j}{h/k^{1/d}}  \right)=B\left(\frac{x-x_0}{2h} \right) = 1 & \ \text{ for } x \in \bigcup_{j=1}^k \mathcal{C}_{h/2k^{1/d}}(m_j) . 
\end{align*}
We note that this result requires a carefully selected relationship between $h$ and $k$, which guarantees that the squared nuisance bumps under $Q_\lambda$ equal the CATE bumps under $P_\lambda$. This also exploits the flat-top bump functions we use, together with a covariate density that only puts mass at these tops, so that the squared terms are constant and no observations occur elsewhere where the bumps are not equal. Without these choices of bump function and covariate density, the expression in \eqref{eq:delta3} would only be bounded by $h^\gamma$, and so $\delta_3$ would only be bounded by $h^{2\gamma}$; in that case,  the $\delta_3$ term would dominate the Hellinger bound in Lemma \ref{lem:hellbound}, and the resulting minimax lower bound would reduce to the oracle rate $n^{-1/(2+d/\gamma)}$, which is uninformative in the low-smoothness regimes we are considering. \\

Now we move to the distance $\delta_1$, which does not end up depending on $h$ and is somewhat easier to handle. For it we have 
\begin{align*}
\delta_1 %&= \max_\ell \sup_\lambda \int_{\mathcal{Z}_\ell} \frac{(p_\lambda - \overline{p})^2}{p_\lambda p_\ell} d\nu \\
&\leq \left( \frac{2^d k}{h^d}  \right) \max_\ell \sup_\lambda \int_{\mathcal{X}_\ell}  \sum_{a,y} \frac{f(x)^2  }{4 p_\lambda(z)}  \left( h/k^{1/d} \right)^{2\beta} \sum_{j=1}^k B \left( \frac{x-m_j}{h/k^{1/d}}  \right)^2 \ dx \\
&\leq \left( \frac{2^d k}{h^d} \right) \left( \frac{2}{\varepsilon} \right)    \left( h/k^{1/d} \right)^{2\beta}  \max_\ell \int_{\mathcal{X}_\ell}  \sum_{j=1}^k B \left( \frac{x-m_j}{h/k^{1/d}}  \right)^2 \ dx  \\
&\leq \left( \frac{2^d \|B \|_2^2}{\varepsilon/2} \right)  \left( h/k^{1/d} \right)^{2\beta}
\end{align*}
where the first line follows by definition, and since $p_\ell \geq (h/2)^d/k$, and $B \left( \frac{x-m_j}{h/k^{1/d}}  \right)=0$  outside of the cube $\mathcal{C}_{h/k^{1/d}}(m_j)$, which implies that 
$$\left\{ \sum_j \lambda_j B \left( \frac{x-m_j}{h/k^{1/d}}  \right) \right\}^2 %= \sum_{j,\ell} \lambda_j \lambda_\ell B \left( \frac{x-m_j}{h/k^{1/d}}  \right)  B \left( \frac{x-m_\ell}{h/k^{1/d}}  \right) 
= \sum_j \lambda_j^2 B \left( \frac{x-m_j}{h/k^{1/d}}  \right)^2, $$ 
the inequality in the second line since $p_\lambda(z)/f(x) \geq \varepsilon$ and $f(x) \leq 2$ as in \eqref{eq:plamfbd} and \eqref{eq:fbd}, and the last inequality since
\begin{align*}
\int_{\mathcal{X}_\ell } \sum_{j=1}^k  B \left( \frac{x-m_j}{h/k^{1/d}}  \right)^2 \ dx &=  \int_{\mathcal{X}_\ell }  B \left( \frac{x-m_\ell}{h/k^{1/d}}  \right)^2  \ dx \leq \frac{h^d}{k} \int B(u)^2 \ du
%\int_{\mathcal{X}_\ell } \sum_{j=1}^k  b_j(x)^m \ dx &=  \int_{\mathcal{X}_\ell }  b_\ell(x)^m \ dx = \int H\{(x-c_\ell) k^{1/d}\}^m \ dx =  \frac{1}{k} \int H(u)^m \ du \lesssim \frac{1}{k} 
\end{align*}
by a change of variables.  \\

For $\delta_2$ we use a mix of the above logic for $\delta_3$ and $\delta_1$. Note that $({q}_\lambda - p_\lambda)^2$ equals
\begin{align*}
 f(x)^2 & \Bigg[  (a-1/2) \left( h/k^{1/d} \right)^\alpha \sum_{j=1}^k \lambda_j B \left( \frac{x-m_j}{h/k^{1/d}}  \right)  \\
& \hspace{.4in} + (2a-1)(2y-1) \left\{ \left( h/k^{1/d} \right)^{2s} \sum_{j=1}^k B \left( \frac{x-m_j}{h/k^{1/d}}  \right)^2  -  \frac{h^\gamma}{4} B\left( \frac{x-x_0}{2h} \right)   \right\} \Bigg]^2  \\
%& \leq 2 f(x)^2 \Bigg[  \frac{1}{4} \left( h/k^{1/d} \right)^{2\alpha} \sum_{j=1}^k B \left( \frac{x-m_j}{h/k^{1/d}}  \right)^2  +  \left\{ \left( h/k^{1/d} \right)^{2s} \sum_{j=1}^k B \left( \frac{x-m_j}{h/k^{1/d}}  \right)^2  -  \frac{h^\gamma}{4} B\left( \frac{x-x_0}{h} \right)  \right\} ^2 \Bigg] \\
& \leq (1/2) f(x)^2  \left( h/k^{1/d} \right)^{2\alpha} \sum_{j=1}^k B \left( \frac{x-m_j}{h/k^{1/d}}  \right)^2  
\end{align*}
where we used the fact that $(a+b)^2 \leq 2(a^2+b^2)$ and $\{ \sum_j \lambda_j B \left( \frac{x-m_j}{h/k^{1/d}}  \right) \}^2 = \sum_j B \left( \frac{x-m_j}{h/k^{1/d}}  \right)^2$, along with the  same logic as above with $\delta_3$ (ensuring the second term in the square equals zero). Now we have 
\begin{align*}
\delta_2 %&= \max_\ell \sup_\lambda \int_{\mathcal{Z}_\ell} \frac{(p_\lambda - \overline{p})^2}{p_\lambda p_\ell} d\nu \\
%&= \left( \frac{2^d k}{h^d}  \right)  \max_\ell \sup_\lambda \int_{\mathcal{X}_\ell}  \sum_{a,y} \frac{f(x)^2  }{4 p_\lambda(z) } \left( h/k^{1/d} \right)^{2\alpha} \sum_{j=1}^k B \left( \frac{x-m_j}{h/k^{1/d}}  \right)^2 \ dx \\
&\leq \left( \frac{2^d k}{h^d} \right) \left( \frac{2}{\varepsilon} \right)   \left( h/k^{1/d} \right)^{2\alpha} \max_\ell \int_{\mathcal{X}_\ell}  \sum_{j=1}^k B \left( \frac{x-m_j}{h/k^{1/d}}  \right)^2 \ dx 
 \leq \left( \frac{2^d \|B \|_2^2}{\varepsilon/2} \right) \left( h/k^{1/d} \right)^{2\alpha}
\end{align*}
using the exact same logic as for $\delta_1$. Plugging these bounds on $(\delta_1,\delta_2,\delta_3)$ into Lemma \ref{lem:hellbound}, together with the fact that $p_j\leq 2(h/2)^d/k$ and $\sum_{j \in S} p_j \leq 2 (h/2)^d$, yields the result. 
\end{proof}

\bigskip

\subsection{Choice of Parameters \& Final Rate} \label{sec:finalrate}

Finally we detail how the parameters $h$ and $k$ can be chosen to ensure the Hellinger distance from Proposition \ref{prop:hellbound} remains bounded, and use the result to finalize the proof of  Theorem \ref{thm:lowerbound}. \\

\begin{proposition} \label{prop:hellbound2}
Let
%$$ h^\gamma =4 k^{-2s/d} \ \text{ and } \  k^{(4s+d+2sd/\gamma)/d} = C^* n^2$$ 
%$$ h^\gamma = 4\left(h/k^{1/d} \right)^{2s} \ \text{ and } \  k = (C^* n^2 )^{d/(4s+d+2sd/\gamma)} $$ 
$$ \frac{h}{k^{1/d}} = \left( \frac{h^\gamma}{4} \right)^{1/2s} = \left( \frac{1}{C^* n^2} \right)^{\frac{1}{d+ 4s+2sd/\gamma}} $$
for $C^* =  2^{2d/\gamma+5}  C (\|B \|_2^2/\varepsilon)^2 $ and $C$ the constant from Proposition \ref{prop:hellbound}. 
Then under the assumptions of Proposition \ref{prop:hellbound} we have
$$  H^2\left(\int P_\lambda^n \ d\varpi(\lambda), \int Q_\lambda^n \ d\varpi(\lambda) \right) \ \leq \ 1  $$ 
and
$h^\gamma =  4 (\sqrt{C^*}  n)^{-1/\left(1 + \frac{d}{2\gamma} + \frac{d}{4s} \right) }$. 
\end{proposition}

\bigskip

\begin{comment}
\begin{proof}
Plugging in the selected values of $h$ and $k$ to the Hellinger bound from Proposition \ref{prop:hellbound} gives
\begin{align*}
C & \left( \frac{ 2^{d+2} \|B \|_2^2}{\epsilon} \right) \left( \frac{n^2 h^d}{k} \right) \left( k^{-4s/d} + k^{-4\alpha/d} \right) \\
&\leq 2C \left( \frac{ 2^{d+2} \|B \|_2^2}{\epsilon} \right) \left( \frac{n^2 h^d}{k} \right)  k^{-4s/d} \\
&=   \frac{C^* n^2}{k^{(4s+d+2sd/\gamma)/d}} = 1
\end{align*}
where the first inequality follows since $\alpha \geq \beta$.
\end{proof}
\end{comment}

The proof of Proposition \ref{prop:hellbound2} follows directly from Proposition \ref{prop:hellbound}, after plugging in the selected values of $h$ and $k$. Importantly, it (together with the alternative construction for the $\beta > \alpha$ case given in Appendix \ref{sec:model2}) also settles the proof of Theorem \ref{thm:lowerbound} via Lemma \ref{lem:minimax}. This follows  since, with the proposed choices of $h$ and $k$, the Hellinger distance is appropriately bounded so that the  term
$( {1 - \sqrt{\alpha(1-\alpha/4)}})/ 2  = (1-\sqrt{3/4})/{2} \approx 0.067 $
in Lemma \ref{lem:minimax} is a constant (greater than 1/20 for example), while the separation in the CATE at $x_0$, which equals $h^\gamma$, is proportional to $n^{-1/(1 + d/2\gamma + d/4s ) }$ under all $P_\lambda$ and $Q_\lambda$. Therefore this separation is indeed the minimax rate in the low smoothness regime where $s < (d/4)/(1+d/2\gamma)$. Note again that, as discussed in Remark \ref{rem:lowsmooth}, when $s > (d/4)/(1+d/2\gamma)$ the rate $n^{-1/(1 + d/2\gamma + d/4s ) }$ is faster than the usual nonparametric regression rate $n^{-1/( 2 + d/\gamma)}$, and so the standard lower bound construction as in Section 2.5 of  \citet{tsybakov2009introduction} indicates that the slower rate $n^{-1/( 2 + d/\gamma)}$ is the tighter lower bound in that regime. \\

\begin{figure}[h!]
\centering
\includegraphics[width=.75\linewidth]{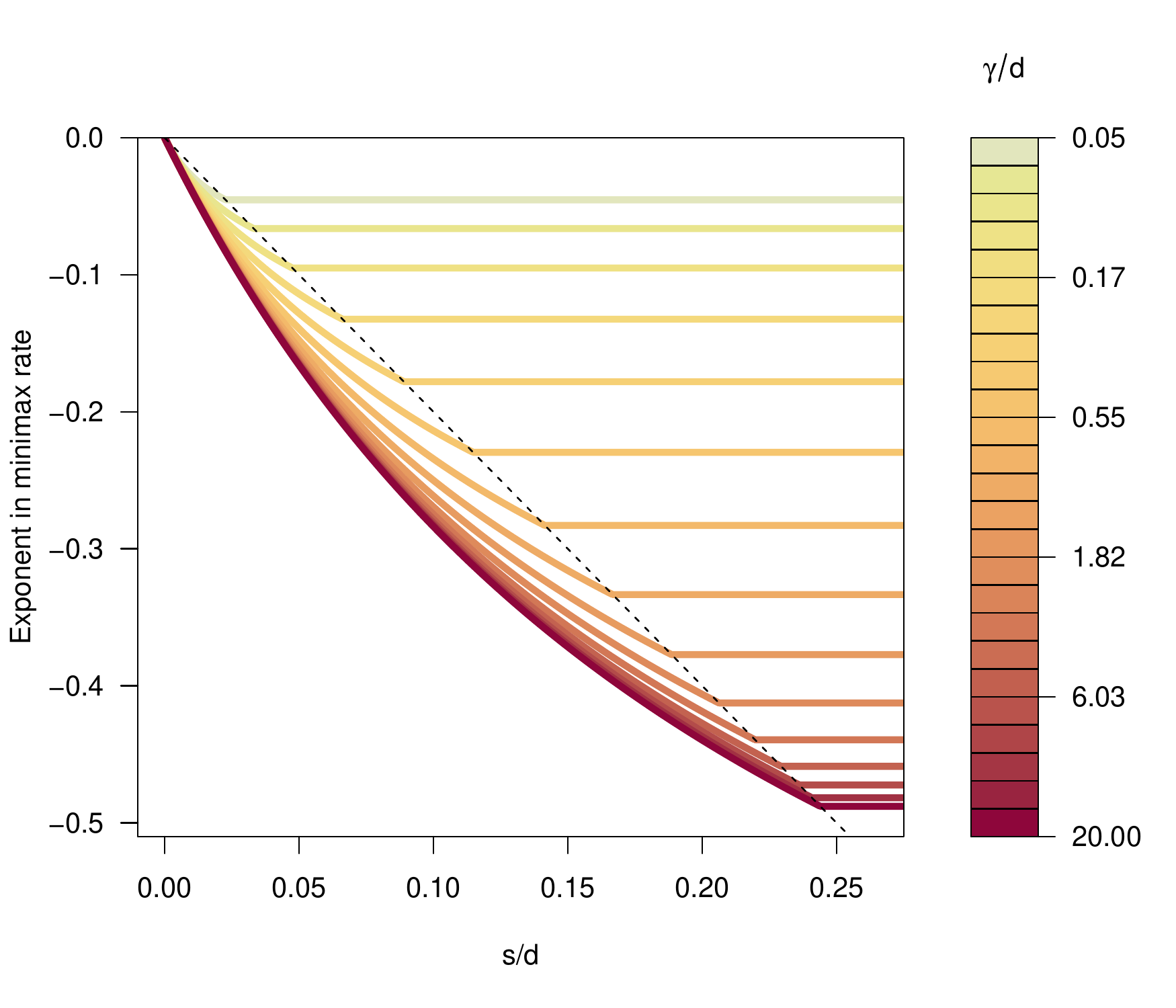}
\caption{The minimax rate for CATE estimation, as a function of average nuisance smoothness $s$ and CATE smoothness $\gamma$, each scaled by covariate dimension $d$. The black dotted line denotes a threshold on the nuisance smoothness $s/d$, below which the oracle nonparametric regression rate $n^{-1/( 2 + d/\gamma)}$ is unachievable (the ``elbow'' phenomenon).}
 \label{fig:rates1}
\end{figure}

Figure \ref{fig:rates1} illustrates the minimax rate from Theorem \ref{thm:lowerbound}, as a function of the average nuisance smoothness $s/d$ (scaled by dimension), and the CATE smoothness scaled by dimension $\gamma/d$.  A number of important features about the rate are worth highlighting. \\

First, of course, the rate never slows with higher nuisance smoothness $s/d$, for any CATE smoothness $\gamma/d$, and vice versa.  However, there is an important elbow phenomenon, akin to that found in functional estimation problems \citep{bickel1988estimating, birge1995estimation, tsybakov2009introduction, robins2009semiparametric}. In particular, the minimax lower bound shows that when the average nuisance smoothness is low enough that $s < \frac{d/4}{1+d/2\gamma}$, the oracle rate $n^{-1/(2+d/\gamma)}$ (which could be achieved if one actually observed the potential outcomes) is in fact unachievable.  This verifies a conjecture in \citet{kennedy2023towards}.  \\

Notably, though, the elbow phenomenon we find in the problem of CATE estimation differs quite substantially from that for classic pathwise differentiable functionals. For the latter, the rate is parametric (i.e., $n^{-1/2}$) above some threshold, and nonparametric ($n^{-1/(1+d/4s)}$) below. In contrast, in our setting the rate matches that of \emph{nonparametric regression} above the threshold, and otherwise is a \emph{combination} of nonparametric regression and functional estimation rates. Thus in this problem there are many elbows, with the threshold depending on the CATE smoothness $\gamma$. In particular, our minimax rate below the threshold, 
$$ n^{-1/\left(1 + \frac{d}{2\gamma} + \frac{d}{4s} \right) }  , $$
is a mixture of the nonparametric regression rate $n^{-1/(1+d/2\gamma)}$ (on the squared scale)  and the classic functional estimation rate $n^{-1/(1+d/4s)}$. This means, for example, that in regimes where the CATE is very smooth, e.g., $\gamma \rightarrow \infty$, the CATE estimation problem begins to resemble that of pathwise-differentiable functional estimation, where the elbow occurs at $s> d/4$, with rates approaching  the parametric rate $n^{-1/2}$ above, and the functional estimation rate $n^{-1/(1+d/4s)}$ below. At the other extreme, where the CATE does not have any extra smoothness, so that $\gamma \rightarrow \beta$ (note we must have $\gamma \geq \beta$), the elbow threshold condition holds for \emph{any} $\alpha \geq 0$. Thus, at this other extreme, there is no elbow phenomenon, and the CATE estimation problem resembles that of smooth nonparametric regression, with optimal rate $n^{-1/(2+d/\beta)}$. For the arguably more realistic setting, where the CATE smoothness $\gamma$ may take intermediate values between $\beta$ and $\infty$, the minimax rate is a mixture, interpolating between the two extremes. All of this quantifies the sense in which the CATE can be viewed as a regression/functional hybrid.  \\

It is also worth mentioning that no estimator previously proposed in the literature (that we know of) attains the minimax rate in Theorem \ref{thm:lowerbound} in full generality. Some estimators have been shown to attain the oracle rate $n^{-1/(2+d/\gamma)}$, but  only under  stronger assumptions than the minimal condition we find here, i.e., that $s > \frac{d/4}{1+d/2\gamma}$. One exception is the undersmoothed R-learner estimator analyzed in \citet{kennedy2023towards}, which did achieve the rate $n^{-1/(2+d/\gamma)}$ whenever $s > (d/4)/(1+d/2\gamma)$, under some conditions (including that $\alpha \geq \beta$). However, in the low-smoothness regime where $s < (d/4)/(1+d/2\gamma)$, that estimator's rate was $n^{-2s/d}$, which is slower than the minimax rate we find here. This motivates our work in the following section, where we propose and analyze a new estimator, whose error matches the minimax rate in much greater generality (under some conditions, e.g., on how well the covariate density is estimated). \\

\begin{remark}
A slightly modified version of our construction also reveals that, when the CATE $\tau(x)=\tau$ is constant, the classic functional estimation rate $n^{-1/\left(1 + \frac{d}{4s}\right)}$ acts as a minimax lower bound. To the best of our knowledge, this result has not been noted elsewhere. \\
\end{remark}

\section{Attainability}
\label{sec:upperbound}

In this section we show that the minimax lower bound of Theorem \ref{thm:lowerbound} is actually attainable, via a new local polynomial version of the R-Learner \citep{nie2017quasi, kennedy2023towards}, based on an adaptation of higher-order influence functions \citep{robins2008higher, robins2009quadratic, robins2017minimax}.  \\

\subsection{Proposed Estimator \& Decomposition}

In this subsection we first describe our proposed estimator, and then give a preliminary error bound, which motivates the specific bias and variance calculations in following subsections. In short, the estimator is a higher-order influence function-based version of the local polynomial R-learner analyzed in \citet{kennedy2023towards}. At its core, the R-Learner essentially regresses outcome residuals on treatment residuals to estimate a weighted average of the CATE. Early versions for a constant or otherwise parametric CATE were studied by \citet{chamberlain1987asymptotic, robinson1988root}, and \citet{robins1994correcting}, with more flexible series, RKHS, and lasso versions studied more recently by \citet{robins2008higher}, \citet{nie2017quasi}, and \citet{chernozhukov2017orthogonal}, respectively. None of this previous work obtained the minimax optimal rates in Theorem \ref{thm:upperbound}.  \\

\begin{definition}[Higher-Order Local Polynomial R-Learner] \label{def:estimator} 
Let $K_h(x) = \frac{1}{h^d} \one(\|x-x_0\| \leq h/2)$. For each covariate $x_j$, $j=1,\dots,d,$ define 
$ \rho(x_j) = \{ \rho_0(x_j), \rho_1(x_j), \dots, \rho_{\lfloor \gamma \rfloor}(x_j) \}^\T $ 
as the first $({\lfloor \gamma \rfloor}+1$) terms of the  Legendre polynomial series (shifted to be orthonormal on  $[0,1]$), 
$$ \rho_{m}(x_j) = \sum_{\ell=0}^m (-1)^{\ell+m} \sqrt{2m+1}  {m \choose \ell} { m + \ell \choose \ell} x_j^\ell . $$
Define $\rho(x)$ to be the corresponding tensor product of all interactions of $\rho(x_1),\dots,\rho(x_d)$ up to order $\lfloor \gamma \rfloor$, which has length  $q = {{ d + \lfloor \gamma \rfloor } \choose { \lfloor\gamma \rfloor }}$ and is orthonormal on $[0,1]^d$, and finally define $\rho_h(x) = \rho\left( 1/2+ (x-x_0)/h \right)$. 
The proposed estimator is then defined as
\begin{equation} \label{eq:estimator}
\widehat\tau(x_0) = \rho_h(x_0)^\T \widehat{Q}^{-1} \widehat{R} 
\end{equation}
where $\widehat{Q}$ is a $q \times q$ matrix and $\widehat{R}$ a $q$-vector given by
\begin{align*}
\widehat{Q} &= \Pn \Big\{  \rho_h(X) K_h(X)\widehat\varphi_{a1}(Z)  \rho_h(X)^\T  \Big\} \\
& \hspace{.4in} + \Un \Big\{ \rho_h(X_1) K_h(X_1) \widehat\varphi_{a2}(Z_1,Z_2) K_h(X_2) \rho_h(X_1)^\T \Big\} \\
\widehat{R} &=  \Pn \Big\{  \rho_h(X_1) K_h(X_1)  \widehat\varphi_{y1}(Z_1) \Big\} + \Un \Big\{  \rho_h(X_1) K_h(X_1) \widehat\varphi_{y2}(Z_1,Z_2) K_h(X_2) \Big\} , 
%\widehat{Q} &= \Un\left[ \rho(X_1) K_h(X_1) \Big\{ \widehat\varphi_{a1}(Z_1) + \widehat\varphi_{a2}(Z_1,Z_2) K_h(X_2) \Big\} \rho(X_1)^\T \right] \\
%\widehat{R} &=  \Un \left[ \rho(X_1) K_h(X_1) \Big\{ \widehat\varphi_{y1}(Z_1) +  \widehat\varphi_{y2}(Z_1,Z_2) K_h(X_2) \Big\} \right]
\end{align*}
respectively, and 
\begin{align*}
\widehat\varphi_{a1}(Z) &= A\{A-\widehat\pi(X)\} \\
\widehat\varphi_{y1}(Z) &= \{Y - \widehat\mu_0(X) \} \{A-\widehat\pi(X)\} \\
\widehat\varphi_{a2}(Z_1,Z_2) &= -\{ A_1 - \widehat\pi(X_1)\} b_{hk}(X_1)^\T \widehat\Omega^{-1} b_{hk}(X_2) A_2  \\
\widehat\varphi_{y2}(Z_1,Z_2) &= -\{A_1 - \widehat\pi(X_1)\} b_{hk}(X_1)^\T \widehat\Omega^{-1} b_{hk}(X_2) \{Y_2 - \widehat\mu_0(X_2)\}  \\
b_{hk}(x) &=  b\{1/2 + (x-x_0)/h\} \one(\|x-x_0 \| \leq h/2) \\
\widehat\Omega &= \int_{v \in [0,1]^d} b(v) b(v)^\T \ d\widehat{F}(x_0 + h(v-1/2)) 
\end{align*}
for $b: \R^d \mapsto \R^k$ a basis of dimension $k$. 
The nuisance estimators $( \widehat{F}, \widehat\pi,\widehat\mu_0)$ are constructed from a separate training sample $D^n$, independent of that on which $\Un$ operates. 
\end{definition}

\bigskip

The estimator in Definition \ref{def:estimator} can be viewed as a \emph{localized} higher-order estimator, and depends on two main tuning parameters: the bandwidth $h$, which controls how locally one averages near $x_0$, and the dimension $k$ of the basis $b$, which controls how bias and variance are balanced in the second-order U-statistic terms in $\widehat{Q}$ and $\widehat{R}$. Specific properties of the basis $b$ are discussed shortly, e.g., just prior to Remark 7 and in \eqref{eq:holder}. We also note that while the basis $b$ will have dimension $k$ growing with sample size, the dimension of the basis $\rho$ is fixed depending on the smoothness of the CATE; for the latter we use the Legendre series for concreteness, but expect other bases to work as well. \\

The U-statistic terms are important for debiasing the first-order sample average terms. In addition, our proposed estimator can be viewed as estimating a locally weighted projection parameter $\tau_h(x_0) = \rho_h(x_0)^\T \theta $, with coefficients given by
\begin{equation} \label{eq:projection}
\argmin_\beta \E \left[ K_h(x) \pi(x)\{1-\pi(x)\} \Big\{ \tau(x) - \beta^\T \rho_h(x) \Big\}^2 \right] = Q^{-1} R
\end{equation}
for 
\begin{align*}
Q &= %\int \rho_h(x) K_h(x) \varphi_{a1}(z) \rho_h(x)^\T \ d\Pb(z) = 
\int \rho_h(x) K_h(x) \pi(x) \{ 1-\pi(x)\} \rho_h(x)^\T \ dF(x)  \\
R &=  %\int \rho_h(x) K_h(x) \varphi_{y1}(z) \ d\Pb(z) =
\int \rho_h(x) K_h(x) \pi(x) \{1-\pi(x)\} \tau(x) \ dF(x) .
\end{align*}
In other words, this projection parameter $\tau_h(x_0)$ is a $K_h(x) \pi(x)\{1-\pi(x)\}$-weighted least squares projection of the CATE $\tau(x)$ on the scaled Legendre polynomials $\rho_h(x)$. Crucially, since $\rho_h(x)$ includes polynomials in $x$ up to order $\lfloor \gamma \rfloor$,  the projection parameter is within $h^\gamma$ of the target CATE; this is formalized in the following proposition. \\

\begin{proposition} \label{prop:projbias}
Let $\tau_h(x) = \rho_h(x)^\T Q^{-1} R$ denote the $x_0$-specific projection parameter from \eqref{eq:projection}, and
assume:
\begin{enumerate}
\item $\tau(x)$ is $\gamma$-smooth,
\item the eigenvalues of $Q$ are bounded below away from zero, and
\item $\int \one\{ \|x-x_0 \| \leq h/2 \} \ dF(x) \lesssim h^d$. 
\end{enumerate}
Then for any $x$ with $\| x - x_0 \| \leq Ch$ we have
$$ | \tau_h(x) - \tau(x) | \lesssim h^\gamma . $$
\end{proposition}

\bigskip

\begin{proof}
This proof follows from a higher-order kernel argument (e.g., Proposition 1.13 of \citet{tsybakov2009introduction}, Proposition 4.1.5 of \citet{gine2021mathematical}), after noting that we can treat $K_h(x)\pi(x)\{1-\pi(x)\}$ itself as a kernel. A similar result was also proved in \citet{kennedy2023towards}. A detailed proof is given in Appendix \ref{app:projbias}. 
\end{proof}

\bigskip

In Proposition \ref{prop:qeigen} in the Appendix we give simple sufficient conditions under which the eigenvalues of $Q$ are bounded. In short, this holds under standard boundedness conditions on the propensity score and covariate density. \\

As mentioned above, our estimator \eqref{eq:estimator} can be viewed as a modified higher-order estimator. More specifically, $\widehat{S}=\widehat{R}-\widehat{Q} \theta$ is closely related to a second-order estimator of the moment condition
\begin{equation} 
\E\Big[ \rho_h(X) K_h(X) \Big\{A-\pi(X) \Big\} \Big\{ Y - \mu_0(X) - A \tau(X)  \Big\} \Big]  = R-Q\theta = 0 \label{eq:moment}
\end{equation}
under the assumption that $\tau(x) = \rho_h(x)^\T \theta$ (this is not exactly true in our case, but  it is enough that it is approximately true locally near $x_0$, as will be proved shortly). Indeed, letting $\widehat{Q}_1=\Pn \{ \rho_h(X) K_h(X)  \widehat\varphi_{a1}(Z) \rho_h(X)^\T \}$ and $\widehat{R}_1=\Pn \{ \rho_h(X) K_h(X)  \widehat\varphi_{y1}(Z) \}$ denote the first terms in $\widehat{Q}=\widehat{Q}_1+\widehat{Q}_2$ and  $\widehat{R}=\widehat{R}_1+\widehat{R}_2$, respectively, we see that $\widehat{S}_1 = \widehat{R}_1 - \widehat{Q}_1 \theta$ is a usual first-order influence function-based estimator of the moment condition \eqref{eq:moment}. Similarly, the second terms $\widehat{Q}_2$ and $\widehat{R}_1$ are akin to the second-order U-statistic corrections that would be added using the higher-order influence function methodology developed by \citet{robins2008higher, robins2009quadratic, robins2017minimax}. However,  these terms differ in two important ways, both relating to localization near $x_0$. First,  the U-statistic is localized with respect to both $X_1$ and $X_2$, i.e., the product $K_h(X_1)K_h(X_2)$ is included, whereas only $K_h(X_1)$ would arise if the goal were purely to estimate the moment condition \eqref{eq:moment} for fixed $h$. Second, the basis functions 
$$ b_{hk}(x) = b\left(1/2 + \frac{x-x_0}{h} \right) \one(\|x-x_0 \| \leq h/2) $$ 
appearing in $\widehat\varphi_{a2}$, $\widehat\varphi_{y2}$, and $\widehat\Omega$ are localized; they only operate on $X$s near $x_0$, stretching them out so as to map the cube $[x_0-h/2,x_0+h/2]^d$ around $x_0$ to the whole space $[0,1]^d$ (e.g., $b_{hk}(x_0-h/2)=b(0)$, $b_{hk}(x_0)=b(1/2)$, etc.). This is the same localization that is used with the Legendre basis $\rho(x)$. In this sense, these localized basis terms spend all their approximation power locally rather than globally away from $x_0$. (Specific approximating properties we require of $b$ will be detailed shortly, in \eqref{eq:holder}). These somewhat subtle distinctions play a crucial role in appropriately controlling bias, as will be described in more detail shortly. \\

\begin{remark}
Note again that, as with other higher-order  estimators, the estimator \eqref{eq:estimator} depends on an initial estimate of the covariate distribution $F$ (near $x_0$), through $\widehat\Omega$. Importantly, we do not take this estimator $\widehat{F}$ to be the empirical distribution, in general, since then our optimal choices of the tuning parameter $k$  would yield $\widehat\Omega$ non-invertible; this occurs with higher-order estimators of pathwise differentiable functionals as well \citep{mukherjee2017semiparametric}. As discussed  in Remark \ref{rem:cond1}, and in more detail shortly, we do give conditions under which the estimation error in $\widehat\Omega$ or $\widehat{F}$ does not impact the overall rate of $\widehat\tau(x_0)$.  \\
\end{remark}

Crucially, Proposition \ref{prop:projbias} allows us to focus on understanding the estimation error in $\widehat\tau(x_0)$ with respect to the projection parameter $\tau_h(x_0)$, treating $h^\gamma$ as a separate approximation bias. The next result gives a finite-sample bound on this error, showing how it is controlled by the error in estimating  the components of $Q$ and $R$.  \\

\begin{proposition} \label{prop:qrerror}
Let $S = R - Q (Q^{-1}R)$ and $\widehat{S} = \widehat{R} - \widehat{Q} (Q^{-1}R)$. The estimator \eqref{eq:estimator} satisfies
$$ \left| \widehat\tau(x_0) - \tau_h(x_0) \right| \ \leq \ \|\rho(1/2) \| \left( \| Q^{-1} \| +   \| \widehat{Q}^{-1} - {Q}^{-1} \| \right) \left\| \widehat{S} - S \right\|   , $$
and if 
 $\| Q^{-1} \|$ and $ \| \widehat{Q}^{-1} - {Q}^{-1} \|$  are bounded above, 
then
\begin{align*}
\E \left| \widehat\tau(x_0) - \tau_h(x_0) \right| \ &\lesssim \  \max_j  \sqrt{ \E\left\{ \E(\widehat{S}_j - S_j \mid D^n)^2 + \var(\widehat{S}_j \mid D^n) \right\} } .
\end{align*}
for $D^n$  a separate independent training sample on which $( \widehat{F}, \widehat\pi,\widehat\mu_0)$ are estimated. 
\end{proposition}

\bigskip

Thus Proposition \ref{prop:qrerror} tells us that bounding the conditional bias and variance of $\widehat{S}=\widehat{R}-\widehat{Q} (Q^{-1}R)$ will also yield finite-sample bounds on the error in $\widehat\tau(x_0)$, relative to the projection parameter $\tau_h(x_0)$. These bias and variance bounds will be derived in the following two subsections. \\

\subsection{Bias}
\label{sec:bias}

In this subsection we derive bounds on  the conditional bias of the estimator $\widehat{S}=\widehat{R}-\widehat{Q} (Q^{-1}R)$, relative to the components of the projection parameter \eqref{eq:projection}, given the training sample $D^n$. The main ideas behind the approach are to use localized versions of higher-order influence function arguments,   along with a specialized localized basis construction, which results in smaller bias due to the fact that the bases only need to be used  in a shrinking window around $x_0$.  \\

Here we rely on the basis $b(x)$ having  
optimal \Holder{} approximation properties, In particular, we assume the approximation error of projections in $L_2$ norm satisfies
\begin{equation} \label{eq:holder}
\| (I-\Pi_b) g \|_{F^*} \ \lesssim \ k^{-s/d}  \ \text{ for any $s$-smooth function $g$}
 \end{equation}
where $\Pi_b g = \argmin_{\ell=\theta^\T b} \int (g - \ell)^2 \ dF^*$ is the usual linear projection of $g$ on $b$, 
for $dF^*(v)=dF(x_0 + h(v-1/2))$ the distribution in  $\mathcal{B}_h(x_0)$, the $h$-ball around $x_0$,  mapped to $[0,1]^d$. 
In a slight abuse to ease notation, we omit the dependence of $\Pi_b g$ on $F^*$.
The approximating condition \eqref{eq:holder} holds for numerous bases, including spline, CDV wavelet, and local polynomial partition series (and polynomial and Fourier series, up to log factors); it is used often in the literature. We refer to \citet{belloni2015some} for more discussion and specific examples (see their Condition A.3 and subsequent discussion in, for example, their Section 3.2). \\

\begin{proposition} \label{prop:bias}
%Let $\| f \|_{K_h}^2 = \int f(x)^2 K_h(x) \ dF(x)$ denote the (squared) kernel-weighted $L_2(\Pb)$ norm. 
Assume:
\begin{enumerate}
\item $\lambda_{\max}(\Omega)$ is bounded above, 
%\item $\|dF^*/dF\|_\infty$ is bounded above, where $F^*(v) = F(x_0 + h(v-1/2))$, 
\item the basis $b$ satisfies approximating condition \eqref{eq:holder}, 
\item $\widehat\pi(x) - \pi(x)$ is $\alpha$-smooth, 
\item $\widehat\mu_0(x) - \mu_0(x)$ is $\beta$-smooth.
\end{enumerate}
Then
\begin{align*}
| \E(\widehat{S}_j - S_j \mid D^n) | & \ \lesssim \ \left( h/k^{1/d} \right)^{2s}  + h^\gamma  \left( h/k^{1/d} \right)^{\alpha} \\
& \hspace{.5in} + \Big( h^\gamma +\| \widehat\mu_0 - \mu_0 \|_{F^*} \Big) \left( \| \widehat\pi - \pi \|_{F^*}   \| \widehat\Omega^{-1}-\Omega^{-1} \|    \right) .
\end{align*}
\end{proposition}

\bigskip

Before delving into the proof, we give some brief discussion. The bias consists of three terms; the first two are the main bias terms that would result even if the covariate distribution $F$ were known, and the third is essentially the contribution from having to estimate $F$. (In Lemma \ref{lem:eigen} of the Appendix we show how the operator norm error of $\widehat\Omega$ is bounded above by estimation error of the distribution $F$ itself). We note that the second of the main bias terms $h^\gamma (h/k^{1/d})^\alpha$ will be of smaller order in all regimes we consider. Compared to the main bias term  in a usual higher-order influence function analysis, which is $k^{-2s/d}$ (e.g., for the average treatment effect), our bias term is smaller;  this is a result of  using the  localized basis $b_{hk}(x)$ defined in \eqref{eq:estimator}, which only has to be utilized locally near $x_0$ (this smaller bias will be partially offset by a larger variance, as discussed in the next subsection). As mentioned in Remark \ref{rem:cond1}, the contribution from having to estimate $F$ is only a third-order term, since the estimation error of $\widehat\Omega$ (in terms of  operator norm) is multiplied by a product  of propensity score errors  (in $L_2(F^*)$ norm) with the sum of regression errors and smoothing bias $h^\gamma$, which is typically of smaller order. In Proposition \ref{prop:bias2}, given after the following proof of Proposition \ref{prop:bias}, we show how the bias simplifies when $F$ is estimated accurately enough.    \\

\begin{proof}
By iterated expectation, the conditional mean of the first-order term in $\widehat{R}$,  i.e., the quantity $\E\{ \rho_h(X_1)  K_h(X_1) \widehat\varphi_{y1}(Z) \mid D^n \} $ is equal to
\begin{align} 
 R& + \int  \rho_h(x) K_h(x)  \{ \pi(x) - \widehat\pi(x) \}\{ \pi(x) \tau(x) + \mu_0(x) - \widehat\mu_0(x) \}  \ dF(x) \nonumber \\
& = R + \int \rho(v) \{ \pi^*(v) - \widehat\pi^*(v) \}\{ \pi^*(v) \tau^*(v) +  \mu_0^*(v) - \widehat\mu_0^*(v) \}  \ dF^*(v) \nonumber
\end{align}
where we use the change of variable $v = \frac{1}{2} + \frac{x-x_0}{h}$ and again define for any function $g: \R^d \mapsto \R$ its corresponding stretched version as $g^*(v) = g(x_0 + h(v-1/2))$. To ease notation it is left implicit that any integral over $v$ is only over $\{v : \|v-1/2\| \leq 1/2 \}$.  Similarly for $\widehat{Q}$ we have that $\E\{ \rho_h(X_1)  K_h(X_1) \widehat\varphi_{a1}(Z) \rho_h(X_1)^\T \mid D^n \} $ equals
\begin{align} 
%\E\{ \rho_h(X_1)  K_h(X_1) & \widehat\varphi_{a1}(Z) \rho_h(X_1)^\T \mid D^n \} 
 Q &+ \int  \rho_h(x) K_h(x)  \{ \pi(x) - \widehat\pi(x) \} \pi(x) \rho_h(x)^\T \ dF(x) \nonumber \\
& = Q + \int \rho(v) \{ \pi^*(v) - \widehat\pi^*(v) \} \pi^*(v) \rho(v) \ dF^*(v) \nonumber
\end{align}
so that for the first-order term in $\widehat{S}$ (denoted $\widehat{R}_1-\widehat{Q}_1 \theta$ in discussion of the moment condition \eqref{eq:moment}) we have
\begin{align}
\E\{ \rho_h(X_1)  &K_h(X_1)  \widehat\varphi_{y1}(Z) \mid D^n \} - \E\{ \rho_h(X_1)  K_h(X_1) \widehat\varphi_{a1}(Z) \rho_h(X_1) \mid D^n \} \theta  \nonumber\\
&= R-Q\theta +  \int \rho(v) \{ \pi^*(v) - \widehat\pi^*(v) \}\{  \mu_0^*(v) - \widehat\mu_0^*(v) \}  \ dF^*(v) \nonumber \\
& \hspace{.5in} + \int \rho(v)\pi^*(v)  \{ \pi^*(v) - \widehat\pi^*(v) \} \{ \tau^*(v) - \tau_h^*(v) \}  \ dF^*(v) . \label{eq:firstbias}
\end{align}
The conditional mean of the second-order influence function term in $\widehat{R}$ is
\begin{align}
&\E\{   \rho_h(X_1) K_h(X_1) \widehat\varphi_{y2}(Z_1,Z_2)K_h(X_2) \mid D^n \} \label{eq:r2ndorder} \\
%&= - \iint  \rho_h(x_1) K_h(x_1) \{\pi(x_1) - \widehat\pi(x_1)\} b_{hk}(x_1)^\T \widehat\Omega^{-1} b_{hk}(x_2) \{ \pi(x_2) \tau(x_2) + \mu_0(x_2) - \widehat\mu_0(x_2)\} K_h(x_2) \ dF(x_2) \ dF(x_1) \nonumber \\
%&= - \iint \rho(v_1)  \{\pi^*(v_1) - \widehat\pi^*(v_1)\} b(v_1)^\T \widehat\Omega^{-1} b(v_2) \{ \pi^*(v) \tau^*(v) + \mu_0^*(v_2) - \widehat\mu_0^*(v_2)\}  \ dF^*(v_2) \ dF^*(v_1) \nonumber \\
&= - \int \rho(v_1)\{\pi^*(v_1) - \widehat\pi^*(v_1)\} \widehat\Pi_{b}(\pi^* \tau^* + \mu_0^*-\widehat\mu_0^*)(v_1) \ dF^*(v_1) \nonumber
\end{align}
where we define
\begin{align*}
\Pi_{b}g^*(u) &= b(u)^\T \Omega^{-1} \int {b}(v) g^*(v) \ dF^*(v) 
\end{align*}
as the $F^*$-weighted linear projection of $g^*$ on the basis $b$, and $\widehat\Pi_{b} g^*(u)$ as the estimated version, which simply replaces $\Omega$ with $\widehat\Omega$. Similarly for $\widehat{Q}\theta$ we have
\begin{align}
&\E\{   \rho_h(X_1) K_h(X_1) \widehat\varphi_{a2}(Z_1,Z_2)K_h(X_2) \rho_h(X_2)^\T \theta \mid D^n \} \label{eq:q2ndorder}  \\
%&= - \iint  \rho_h(x_1) K_h(x_1) \{ \pi(x_1) - \widehat\pi(x_1)\} b_{hk}(x_1)^\T \widehat\Omega^{-1} b_{hk}(x_2)  \pi(x_2)  K_h(x_2) \rho_h(x_2)^\T \theta \ dF(x_2) \ dF(x_1) \nonumber \\
%&= - \iint \rho(v_1)  \{ \pi^*(v_1) - \widehat\pi^*(v) \} b(v_1)^\T \widehat\Omega^{-1} b(v_2) \pi^*(v_2) \rho(v_2)^\T  \theta \ dF^*(v_2) \ dF^*(v_1) \nonumber \\
&= - \int \rho(v_1) \{ \pi^*(v_1) - \widehat\pi^*(v_1) \} \widehat\Pi_{b}(\pi^* \tau_h^*)(v_1) \ dF^*(v_1)  \nonumber
\end{align}
so that the conditional mean \eqref{eq:r2ndorder} minus the  conditional mean \eqref{eq:q2ndorder} equals
\begin{align}
 %\E\{  & \rho_h(X_1) K_h(X_1) \widehat\varphi_{y2}(Z_1,Z_2)K_h(X_2) \mid D^n \} - \E\{   \rho_h(X_1) K_h(X_1) \widehat\varphi_{a2}(Z_1,Z_2)K_h(X_2) \rho_h(X_2)^\T \theta \mid D^n \} \nonumber \\
 - \int \rho(v_1)\{\pi^*(v_1) - \widehat\pi^*(v_1)\} \widehat\Pi_{b}\{ \pi^* ( \tau^* - \tau_h^* ) + (\mu_0^*-\widehat\mu_0^*)\}(v_1) \ dF^*(v_1) . \label{eq:secondbias}
\end{align}
Therefore adding the first- and second-order expected values in \eqref{eq:firstbias} and \eqref{eq:secondbias},  the overall bias relative to $S$ is
\begin{align}
&\int  \rho(v)  \{ \pi^*(v) - \widehat\pi^*(v)\}  (I-\widehat\Pi_{b}) \{ \pi^* ( \tau^* - \tau_h^* ) + (\mu_0^*-\widehat\mu_0^*)\}(v) \ dF^*(v) \nonumber \\
%&= \int \rho(v) \{ \pi^*(v) - \widehat\pi^*(v)\} (I-\Pi_{b}) \{ \pi^* ( \tau^* - \tau_h^* ) + (\mu_0^*-\widehat\mu_0^*)\}(v) \ dF^*(v) \nonumber \\
%& \hspace{.6in} + \int \rho(v) \{ \pi^*(v) - \widehat\pi^*(v)\} (\Pi_{b}-\widehat\Pi_{b}) \{ \pi^* ( \tau^* - \tau_h^* ) + (\mu_0^*-\widehat\mu_0^*)\}(v) \ dF^*(v)  \nonumber \\
&= \int (I-\Pi_b) \{ \rho (\pi^* - \widehat\pi^*)\}(v) (I-\Pi_{b})\{ \pi^* ( \tau^* - \tau_h^* ) + (\mu_0^*-\widehat\mu_0^*)\}(v) \ dF^*(v) \label{eq:mainbias}  \\
& \hspace{.2in} + \int \rho(v) \{ \pi^*(v) - \widehat\pi^*(v)\} (\Pi_{b}-\widehat\Pi_{b})\{ \pi^* ( \tau^* - \tau_h^* ) + (\mu_0^*-\widehat\mu_0^*)\}(v) \ dF^*(v)  \label{eq:covbias} 
\end{align}
by the orthogonality of a projection with its residuals (Lemma \ref{lem:projections}(i)). \\

Now we analyze the bias terms \eqref{eq:mainbias} and \eqref{eq:covbias} separately; the first is the main bias term, which would arise even if the covariate density were known, and the second is the contribution coming from having to estimate the covariate density. \\

Crucially, by virtue of using the localized basis $b_{hk}$, the projections in these bias terms are of stretched versions of the nuisance functions $(\pi^*-\widehat\pi^*)$ and $(\mu_0^*-\widehat\mu_0^*)$, on the standard  \emph{non-localized} basis $b$, with weights equal to the stretched density $dF^*$. This is important because stretching a function increases its smoothness; in particular, the stretched and scaled function $g^*(v)/h^\alpha$ is $\alpha$-smooth whenever $g$ is $\alpha$-smooth. This follows since  $| D^{\lfloor \alpha \rfloor} g^*(v) - D^{\lfloor \alpha \rfloor} g^*(v') |$ equals
\begin{align*}
%\Big| D^{\lfloor \alpha \rfloor} g^*(v) &- D^{\lfloor \alpha \rfloor} g^*(v') \Big| = 
&\left| D^{\lfloor \alpha \rfloor} g( x_0 + h(v-1/2) ) - D^{\lfloor \alpha \rfloor} g( x_0 + h(v'-1/2) ) \right|  \\
&\hspace{.4in}= h^{\lfloor \alpha \rfloor} \left| g^{(\lfloor \alpha \rfloor)}( x_0 + h(v-1/2) )  - g^{(\lfloor \alpha \rfloor)}( x_0 + h(v'-1/2) )  \right| \\
&\hspace{.4in}\lesssim h^\alpha | v - v'| 
\end{align*}
where the second equality follows by the chain rule, and the third since $g$ is $\alpha$-smooth. Thus the above  implies 
$ h^{-\alpha} \left| D^{\lfloor \alpha \rfloor} g^*(v) - D^{\lfloor \alpha \rfloor} g^*(v') \right| \lesssim  | v - v'| ,  $
i.e., that $g^*(v)/h^\alpha$ is $\alpha$-smooth. \\

Therefore if $g$ is $\alpha$-smooth, then 
$\| (I - \Pi_b) g^*/h^\alpha \|_{F^*} \lesssim k^{-\alpha/d}$
by the \Holder{} approximation properties \eqref{eq:holder} of the basis $b$, and so it follows that
\begin{equation}
\| (I - \Pi_b) g^* \|_{F^*} \ \lesssim \ h^\alpha k^{-\alpha/d} = \left(h/k^{1/d} \right)^{\alpha} \label{eq:smallbias}
\end{equation}
for any $\alpha$-smooth function $g$. \\

Therefore now consider the bias term \eqref{eq:mainbias}. 
This term satisfies
\begin{align*}
\int \Big[ &(I- \Pi_{b})\{\rho(\pi^*-\widehat\pi^*)\}(v) \Big] \Big[ (I-\Pi_{b}) \{ \pi^* ( \tau^* - \tau_h^* ) + (\mu_0^*-\widehat\mu_0^*)\}(v) \Big] \ dF^*(v) \\
&\leq \| (I-\Pi_{b}) \{\rho(\pi^*-\widehat\pi^*)\}  \|_{F^*} \Big[ \| (I-\Pi_{b})\{\pi^*(\tau^*-\tau_h^*) \} \|_{F^*}  \\
& \hspace{.4in} + \| (I-\Pi_{b})(\mu_0^*-\widehat\mu_0^*) \|_{F^*}  \Big] \\
&  \lesssim  \  \left(h/k^{1/d}\right)^{\alpha} \Big\{ h^\gamma  + \left(h/k^{1/d}\right)^\beta \Big\} = \left(h/k^{1/d}\right)^{2s} +h^\gamma \left(h/k^{1/d}\right)^\alpha
\end{align*}
where the second line follows by Cauchy-Schwarz, and the third by \eqref{eq:smallbias}, 
since $(\pi-\widehat\pi)$ and  $(\mu_0-\widehat\mu_0)$ are assumed $\alpha$- and $\beta$-smooth, respectively (note $\rho(v)$ is a polynomial, so the smoothness of $\rho(\pi^*-\widehat\pi^*)$ is the same as $(\pi^*-\widehat\pi^*)$), along with the fact that
\begin{align*}
\| (I-\Pi_{b}) \pi^* ( \tau^* &- \tau_h^* ) \|_{F^*}^2 \leq \|  \pi^* ( \tau^* - \tau_h^* ) \|_{F^*}^2 \\
&\leq \int K_h(x) \Big\{ \tau(x) - \tau_h(x) \Big\}^2 dF(x) \lesssim h^{2\gamma}
\end{align*}
where the first inequality follows by Lemma \ref{lem:projections}(ii), the second by definition of $F^*$ and since $\pi(x) \leq 1$, and the last by Proposition \ref{prop:projbias}. \\

Now for the term in  \eqref{eq:covbias}, let $\theta_{b,g} = \Omega^{-1} \int b g \ dF^*$ denote the coefficients of the projection $\Pi_{b}g$, and note for any functions $g_1,g_2$ we have
\begin{align*}
\int g_1(\Pi_{b} - \widehat\Pi_{b})(g_2) \ dF^* %\\ &=  \left\{ \int {b}_h(x) K_h(x) g(x)   \ dF(x) \right\}^\T   (\Omega^{-1} - \widehat\Omega^{-1})  \left\{   \int {b}_h(t) K_h(t) f(t)  \ d\Pb(t) \right\} \\
&= \left( \Omega^{1/2} \theta_{b, g_1} \right)^\T \Omega^{1/2} (\Omega^{-1} - \widehat\Omega^{-1}) \Omega^{1/2}  \left( \Omega^{1/2} \theta_{b, g_2} \right) \\
&\leq \| g_1 \|_{F^*}  \| \Omega^{1/2} (\Omega^{-1} - \widehat\Omega^{-1}) \Omega^{1/2}  \| \| g_2 \|_{F^*} \\
%&= \| f \|_{K_h} \| g \|_{K_h}  \| \Omega^{1/2} \Omega^{-1}  ( \widehat\Omega - \Omega)  \widehat\Omega^{-1} \Omega^{1/2}  \| \\
&\leq  \| g_1 \|_{F^*} \| g_2 \|_{F^*}  \|  \Omega \| \| \widehat\Omega^{-1} - \Omega^{-1} \|  
\end{align*}
where the first equality follows by definition, the second line since the $L_2$ norm of the coefficients of a (weighted) projection is no more than the weighted $L_2(\Pb)$ norm of the function itself (Lemma \ref{lem:projections}(iii)), and the last by the sub-multiplicative property of the operator norm, along with the fact that $\| \Omega^{1/2} \|^2 = \| \Omega \|$. 
\end{proof}

\bigskip

Several of our results require the eigenvalues of $\Omega$ to be bounded above and below away from zero. Proposition \ref{prop:oeigen} in the Appendix gives simple sufficient conditions for this to hold (similar to Proposition \ref{prop:qeigen} for the matrix $Q$, which was mentioned earlier after Proposition \ref{prop:projbias}). \\

The next result  is a refined version of Proposition \ref{prop:bias}, giving high-level conditions under which estimation of $F$ itself (rather than the matrix $\Omega^{-1}$)  does not impact the bias. We refer to Remark \ref{rem:cond1} for more detailed discussion of these conditions, and note that the result follows from Proposition \ref{prop:bias} together with Lemma \ref{lem:eigen} in the Appendix. \\

\begin{proposition} \label{prop:bias2}
Under the assumptions of Proposition \ref{prop:bias2}, and if
\begin{enumerate}
\item $\lambda_{\min}(\Omega)$ is bounded below away from zero, 
\item $\| d\widehat{F}^*/dF^*\|_\infty$ is bounded above and below away from zero, 
\item $\| (d\widehat{F}^*/dF^*)-1 \|_\infty \lesssim \frac{ (h/k^{1/d})^{2s} }{ \| \widehat\pi-\pi \|_{F^*} (\| \widehat\mu_0-\mu_0 \|_{F^*} + h^\gamma) }$,
\end{enumerate}
then, when $h^\gamma \lesssim (h/k^{1/d})^\beta$, the bias satisfies $| \E(\widehat{S}_j - S_j \mid D^n) | \lesssim 
 \left( h/k^{1/d} \right)^{2s}.$
\end{proposition}

\bigskip

\subsection{Variance}
\label{sec:variance}

In this subsection we derive bounds on  the conditional variance of the estimators $\widehat{R}_j$ and $\widehat{Q}_{j\ell}$, given the training sample $D^n$. The main tool used here is a  localized version of second-order U-statistic variance arguments, recognizing that our higher-order estimator is, conditionally, a second-order U-statistic over $nh^d$ observations. \\

\begin{proposition} \label{prop:variance}
Assume:
\begin{enumerate}
\item $y^2$, $\widehat\pi^2$, $\widehat\mu_0^2$, and $\| \widehat\mu_0 - \mu_0 \|_{F^*}$ are all bounded above, and
\item $\lambda_{\max}(\Omega)$ is bounded above. 
\end{enumerate}
Then
\begin{align*}
\var(\widehat{S}_j   \mid D^n)  & \ \lesssim \ \frac{1}{nh^d} \left( 1 + \frac{k  }{n  h^{d}} \left( 1 + \| \widehat\Omega^{-1}-\Omega^{-1} \|^2 \right) \right) .
\end{align*}
\end{proposition}

\bigskip

We give the proof of Proposition \ref{prop:variance} in Appendix \ref{app:propvariance}, and so just make some comments here.  First, the variance here is analogous to that of a  higher-order (quadratic) influence function estimator (cf.\ Theorem 1 of \citet{robins2009quadratic}), except with sample size $n$ deflated to $nh^d$.  This is to be expected given the double localization in our proposed estimator. Another important note is that the contribution to the variance from having to estimate $F$ is relatively minimal, compared to the bias, as detailed in Proposition \ref{prop:bias}. For the bias, non-trivial rate conditions are needed to ensure estimation of $F$ does not play a role, whereas for the variance one only needs the operator norm of $\widehat\Omega^{-1} - \Omega^{-1}$ to be bounded (under regularity conditions, this amounts to the estimator $\widehat{F}$ only having bounded errors, in a relative sense, as noted in the following remark). \\

\begin{remark}
By Lemma \ref{lem:eigen}, under the assumptions of Proposition \ref{prop:bias2}, it follows that 
$$ \| \widehat\Omega^{-1}-\Omega^{-1} \| \ \lesssim \ \| (d\widehat{F}^*/dF^*) - 1 \|_\infty , $$
so estimation of $F$ will not affect the conditional variances as long as the error of $\widehat{F}$ is bounded in uniform norm.  \\
\end{remark}

\subsection{Overall Rate} \label{sec:overallrate}

Combining the approximation bias in Proposition \ref{prop:projbias} with the decomposition in Proposition \ref{prop:qrerror}, and the bias and variance bounds from Proposition \ref{prop:bias2} and Proposition \ref{prop:variance}, respectively, shows that
\begin{equation} \label{eq:overallrate}
 \E_P | \widehat\tau(x_0) - \tau_P(x_0) |  \ \lesssim \ h^\gamma + \left( h/k^{1/d} \right)^{2s} + \sqrt{\frac{1}{nh^d} \left( 1 + \frac{k}{nh^d} \right) } 
 \end{equation}
under all the combined assumptions of these results, which are compiled in the statement of Theorem \ref{thm:upperbound} below. \\

 The first two terms in \eqref{eq:overallrate} are the bias, with $h^\gamma$ an oracle bias that would remain even if one had direct access to the potential outcomes $(Y^1-Y^0)$ (or equivalently,  samples of $\tau(X)+\epsilon$ for some $\epsilon$ with conditional mean zero), and $(h/k^{1/d})^{2s}$ analogous to a squared nuisance bias term, but shrunken due to the stretching induced by the localized basis $b_{hk}$. Similarly, $1/(nh^d)$ is an oracle variance that would remain even if given access to the potential outcomes, whereas the $k/(nh^d)$ factor is a contribution from nuisance estimation (akin to the variance of a series regression on $k$ basis terms with $nh^d$ samples). \\

Balancing bias and variance in \eqref{eq:overallrate} by taking the tuning parameters to satisfy 
$$h \sim n^{-(1/\gamma)/\left(1 + \frac{d}{2\gamma} + \frac{d}{4s} \right) } $$
and
$$ k \sim n^{\left( \frac{d}{2s} - \frac{d}{\gamma} \right)/\left(1 + \frac{d}{2\gamma} + \frac{d}{4s} \right) } $$
ensures the rate matches the minimax lower bound from Theorem \ref{thm:lowerbound} (in the low smoothness regime), proving that lower bound is in fact tight. (In the high smoothness regime where $s>(d/4)/(1+d/2\gamma)$, one can simply take $k \sim nh^d$ and $h \sim n^{-1/(2\gamma+d)}$). This is formalized in the following theorem. \\

\begin{theorem} \label{thm:upperbound}
Assume the regularity conditions:
\begin{enumerate}[label=\Alph*]
\item[A.] The eigenvalues of $Q$ and $\Omega$ are bounded above and below away from zero. \label{cnd:reg1}
\item[B.] $\widehat\pi(x)-\pi(x)$ is $\alpha$-smooth and $\widehat\mu_0(x) - \mu_0(x)$ is $\beta$-smooth. \label{cnd:reg2}
\item[C.] The quantities $y^2$, $(\widehat\pi^2 , \widehat\mu_0^2)$, $\| \widehat\mu_0 - \mu_0 \|_{F^*}$, and $ \| \widehat{Q}^{-1} - {Q}^{-1} \|$  are all bounded above, and $\| d\widehat{F}^*/dF^* \|_\infty$ is bounded above and below away from zero. \label{cnd:reg3}
\end{enumerate}
Also assume the basis $b$  satisfies \Holder{} approximating condition \eqref{eq:holder},  and:
\begin{enumerate}
\item $\| (d\widehat{F}^*/dF^*)-1 \|_\infty \lesssim \frac{n^{-1/\left(1 + \frac{d}{2\gamma} + \frac{d}{4s} \vee \left(1 + \frac{d}{2\gamma}\right) \right) }  }{ \| \widehat\pi-\pi \|_{F^*} (\| \widehat\mu_0-\mu_0 \|_{F^*} + h^\gamma)}$, \label{cnd:covest}
\item $\pi(x)$ is $\alpha$-smooth, and $\epsilon \leq \pi(x) \leq 1-\epsilon$ for some $\epsilon>0$, 
\item $\mu_0(x)$ is $\beta$-smooth, 
\item $\tau(x)$ is $\gamma$-smooth.
\end{enumerate}
Finally let the tuning parameters satisfy
$$  h \sim n^{-(1/\gamma)/\left(1 + \frac{d}{2\gamma} + \frac{d}{4s} \right) } \ \text{ and } \ k \sim n^{\left( \frac{d}{2s} - \frac{d}{\gamma} \right)/\left(1 + \frac{d}{2\gamma} + \frac{d}{4s} \right) }  $$
if $s < \frac{d/4}{1+d/2\gamma}$, or $h \sim n^{-\frac{1}{2\gamma+d} }$ and $k \sim nh^d$ otherwise. 
Then the estimator $\widehat\tau$ from Definition \ref{def:estimator} has error upper bounded as
\begin{align*}
\E_P | \widehat\tau(x_0) - \tau_P(x_0)| \lesssim \begin{cases}
n^{-1/\left(1 + \frac{d}{2\gamma} + \frac{d}{4s} \right) } & \text{ if } s < \frac{d/4}{1+d/2\gamma}  \\
n^{-1/\left(2 + \frac{d}{\gamma}  \right) } & \text{ otherwise.} 
\end{cases}
\end{align*}
\end{theorem}

\bigskip

We refer to Section \ref{sec:finalrate} for more detailed discussion and visualization of the rate from Theorem \ref{thm:upperbound}. Here we give two remarks discussing the regularity conditions A--C and Condition \ref{cnd:covest} (which ensures the covariate distribution is estimated accurately enough).  \\

\begin{remark}
Condition \ref{cnd:reg1} is a standard collinearity restriction used with least squares estimators; simple sufficient conditions are given in Propositions \ref{prop:qeigen} and \ref{prop:oeigen} in the Appendix. In Lemma \ref{lem:covdensity} in the Appendix we also prove that this condition holds for a class of densities contained in the model $\mathcal{P}$ in Theorem \ref{thm:lowerbound}, ensuring that the upper bound  holds over the same submodel. A sufficient condition for Condition \ref{cnd:reg2} to hold is that the estimators $\widehat\pi(x)$ and $\widehat\mu_0(x)$ match the (known) smoothness of $\pi(x)$ and $\mu_0(x)$; this would be the case for standard minimax optimal estimators based on series or local polynomial methods. Condition \ref{cnd:reg3} is a mild boundedness condition on the outcome $Y$ (which could be weakened at the expense of adding some complexity in the analysis), as well as the nuisance estimators $(\widehat{F}^*, \widehat\pi,\widehat\mu_0$), and even weaker, the errors $\| \widehat\mu_0 - \mu_0 \|_{F^*}$ and $ \| \widehat{Q}^{-1} - {Q}^{-1} \|$ (which would typically not only be bounded but decreasing to zero). \\
\end{remark}

\begin{remark} \label{rem:cond1}
First, Condition \ref{cnd:covest} of Theorem \ref{thm:upperbound} will of course hold if the covariate distribution is estimated at a rate faster than that of the CATE (i.e., the numerator of the rate in Condition \ref{cnd:covest}); however, it also holds under substantially weaker conditions, depending on how accurately $\pi$ and $\mu_0$ are estimated. This is because the condition really amounts to a third-order term (the covariate distribution error {multiplied} by a \emph{product} of nuisance errors) being of smaller order than the CATE rate.  Specifically, the result of Theorem \ref{thm:upperbound} can also be written as
\begin{equation} \label{eq:ub2}
 \E_P | \widehat\tau(x_0) - \tau_P(x_0)| \ \lesssim \ n^{-1/\left(1 + \frac{d}{2\gamma} + \frac{d}{4s} \vee \left(1 + \frac{d}{2\gamma}\right) \right) } + R_{3,n} , 
\end{equation}
for the third-order error term 
$$ R_{3,n} = \| (d\widehat{F}^*/dF^*)-1 \|_\infty \| \widehat\pi-\pi \|_{F^*} \Big(\| \widehat\mu_0-\mu_0 \|_{F^*} + h^\gamma \Big) , $$
so that Condition 1 simply requires this third-order term to be of smaller order than the first minimax optimal rate in \eqref{eq:ub2} (note in the above that $h^\gamma$ matches the overall CATE estimation error, under our tuning parameter choices, which would typically be of smaller order than the regression error $\| \widehat\mu_0 - \mu_0 \|_{F^*}$). 
Second, we note that we leave the condition in terms of the $L_2(F^*)$ errors $\| \widehat\pi-\pi \|_{F^*}$ and $\| \widehat\mu_0-\mu_0 \|_{F^*}$ because, although we assume $\pi$ and $\mu_0$ are $\alpha$- and $\beta$-smooth, technically, they do not need to be estimated at particular rates for any of the other results we prove to hold. Of course, under these smoothness assumptions, there are available minimax optimal estimators for which 
$$ \| \widehat\pi-\pi \|_{F^*} \asymp n^{-1/(2+d/\alpha)} \ \text{ and } \ \| \widehat\mu_0-\mu_0 \|_{F^*} \asymp n^{-1/(2+d/\beta)} . $$
If in addition there exists some $\zeta$ for which $\| (d\widehat{F}^*/dF^*)-1 \|_\infty \asymp n^{-1/(2+d/\zeta)}$ (e.g., if $F$ has a density that is $\zeta$-smooth), then Condition 1 reduces to  $ \zeta > {d}/{(1/M_{\alpha,\beta,\gamma,d}-2)}$, for
$$  M_{\alpha,\beta,\gamma,d} \equiv  \frac{1}{1+d/2\gamma + d/4s \vee (1 + d/2\gamma)}  -\frac{1}{2+d/\alpha} -  \frac{1}{2+d/\beta}  . $$
Exploring CATE estimation under weaker conditions on the covariate distribution is an interesting avenue for future work; we suspect the minimax rate changes depending on what is assumed about this distribution, as is the case for average effects (e.g.,  page 338 of \citet{robins2008higher}) and conditional variance estimation \citep{wang2008effect, shen2020optimal}. \\
\end{remark}

\section{Discussion}

In this paper we have characterized the minimax rate for estimating heterogeneous causal effects in a smooth nonparametric model. We derived a lower bound on the minimax rate using a localized version of the method of fuzzy hypotheses, and a matching upper bound via a new  local polynomial R-Learner estimator based on higher-order influence functions. We also characterize how the minimax rate changes depending on whether the propensity score or regression function is smoother, either when one parametrizes the control or the marginal regression function.  The minimax rate has several important features. First, it exhibits a so-called elbow phenomenon: when the nuisance functions (regression and propensity scores) are smooth enough, the rate matches that of standard smooth nonparametric regression (the same that would be obtained if potential outcomes were actually observed). On the other hand, when the average nuisance smoothness is below the relevant threshold, the rate obtained is slower. This leads to a second important feature: in the latter low-smoothness regime, the minimax rate is a mixture of the minimax rates for nonparametric regression and functional estimation. This quantifies how the CATE can be viewed as a regression/functional hybrid. \\

There are numerous important avenues left for future work. We detail a few briefly here, based on: different error metrics, inference and testing, adaptivity, and practical implementation. First, we have focused on estimation error at a point, but one could also consider global rates in $L_2$ or $L_\infty$ norm, for example. We expect $L_2$ rates to be the same, and $L_\infty$ rates to be the same up to log factors, but verifying this would be useful. In addition, it would also be very important to study the distribution of the proposed estimator, beyond just bias and variance, e.g., for purposes of inference ($L_\infty$ rates could also be useful in this respect). Relatedly, one could consider minimax rates for testing whether the CATE is zero, for example, versus $\epsilon$-separated in some distance.  The goal of the present work is mostly to further our theoretical understanding of the fundamental limits of CATE estimation, so there remains lots  to do to make the optimal rates obtained here achievable in practice. For example, although we have specified particular values of the tuning parameters $h$ and $k$ to confirm attainability of our minimax lower bound, it would be practically useful to have more data-driven approaches for selection. In particular, the optimal tuning values depend on underlying smoothness, and since in practice this is often unknown, a natural next step is to study adaptivity. For example one could study whether approaches based on Lepski's method could be used, as in \citet{mukherjee2015lepski} and \citet{liu2021adaptive}. There are also potential computational challenges associated with constructing the tensor products  in $\rho(x)$ when dimension $d$ is not small, as well as evaluating the U-statistic terms of our estimator, and inverting the matrices $\widehat{Q}$ and $\widehat\Omega$. Finally, in this work we have assumed the nuisance functions are \Holder{}-smooth, a classic infinite-dimensional function class from which important insights can be drawn.  However, it will be important to explore minimax rates in other function classes as well.  \\

\section*{Acknowledgements}

EK gratefully acknowledges support from NSF DMS Grant 1810979, NSF CAREER Award 2047444, and  NIH R01 Grant LM013361-01A1, and SB and LW from NSF Grant DMS1713003. EK also thanks Matteo Bonvini and Tiger Zeng for very helpful discussions. \\

\section*{References}
\vspace{-1cm}
\bibliographystyle{abbrvnat}
\bibliography{/Users/kennedye/Desktop/FLASHDRIVE/research/bibliography}

\pagebreak

\appendix

\appendixpage

\section*{Roadmap of Main Results}

The following figures illustrate how Theorems \ref{thm:lowerbound} and \ref{thm:upperbound} follow from supporting propositions and lemmas (with arrows indicating that a result plays a role in implying another). \\

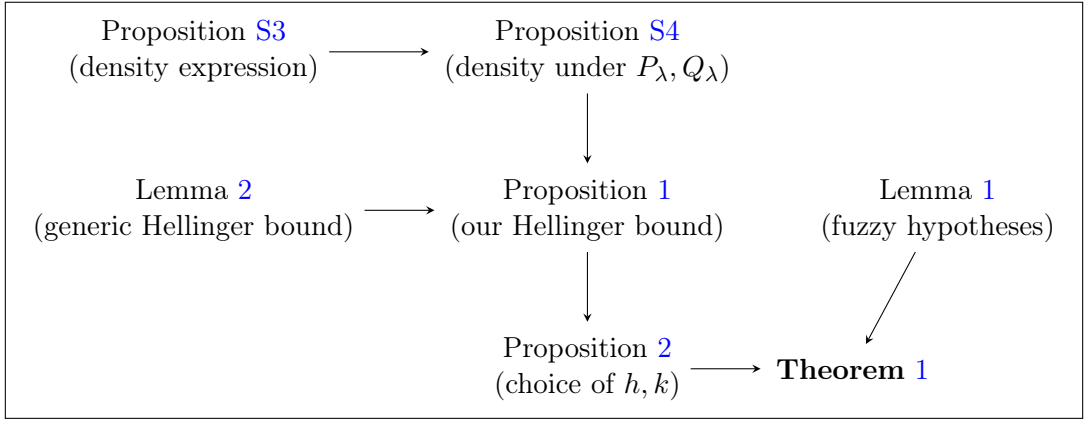
\begin{figure}[h!]
\begin{center}
\fbox{
\begin{tikzpicture}[->, shorten >=2pt,>=stealth, node distance=1cm, align=center]
\node[] (1) {Lemma \ref{lem:hellbound} \\ (generic Hellinger bound)};
\node (3) [above= of 1] {Proposition \ref{prop:dens} \\ (density expression)};
\node (4) [right= of 1] {Proposition \ref{prop:hellbound} \\ (our Hellinger bound)};
\node (2) [above= of 4] {Proposition \ref{prop:dens2} \\ (density under $P_\lambda,Q_\lambda$)};
\path (1) edge node {} (4);
\path (2) edge node {} (4);
\path (3) edge node {} (2);
\node (5) [below= of 4] {Proposition \ref{prop:hellbound2} \\ (choice of $h,k$)};
\path (4) edge node {} (5);
\node (6) [right= of 4] {Lemma \ref{lem:minimax} \\ (fuzzy hypotheses)};
\node (7) [right= of 5] {\textbf{Theorem} \ref{thm:lowerbound}};
\path (5) edge node {} (7);
\path (6) edge node {} (7);
\end{tikzpicture} }
\caption{Roadmap for lower bound result in Theorem \ref{thm:lowerbound}.}
\end{center}
\end{figure}

\bigskip

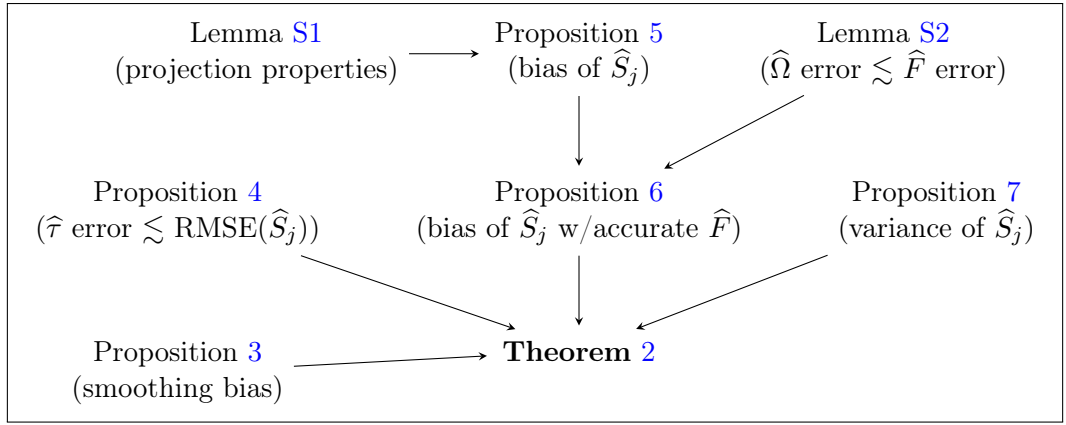
\begin{figure}[h!]
\begin{center}
\fbox{
\begin{tikzpicture}[->, shorten >=2pt,>=stealth, node distance=1cm, align=center]
\node[] (1) {Lemma \ref{lem:projections} \\ (projection properties)};
\node (2) [right= of 1] {Proposition \ref{prop:bias} \\ (bias of $\widehat{S}_j$)};
\node (3) [right= of 2] {Lemma \ref{lem:eigen} \\ ($\widehat\Omega$ error $\lesssim$ $\widehat{F}$ error)};
\node (4) [below= of 2] {Proposition \ref{prop:bias2} \\ (bias of $\widehat{S}_j$ w/accurate $\widehat{F}$)};
\path (1) edge node {} (2);
\path (2) edge node {} (4);
\path (3) edge node {} (4);
\node (5) [right= of 4] {Proposition \ref{prop:variance} \\ (variance of $\widehat{S}_j$)};
\node (6) [left= of 4] {Proposition \ref{prop:qrerror} \\ ($\widehat\tau$ error $\lesssim$ RMSE($\widehat{S}_j$))};
\node (7) [below= of 6] {Proposition \ref{prop:projbias} \\ (smoothing bias)};
\node (8) [below= of 4] {\textbf{Theorem} \ref{thm:upperbound}};
\path (4) edge node {} (8);
\path (5) edge node {} (8);
\path (6) edge node {} (8);
\path (7) edge node {} (8);
\end{tikzpicture} }
\caption{Roadmap for upper bound result in Theorem \ref{thm:upperbound}.}
\end{center}
\end{figure}

\bigskip

Propositions \ref{prop:qrerror2} and \ref{prop:bias_alt} are not included in the above figures since they support the alternative parametrization detailed in Theorem \ref{thm:lowerbound2}; however they play  roles analogous to Propositions \ref{prop:qrerror} and \ref{prop:bias}, respectively. Propositions \ref{prop:qeigen} and \ref{prop:oeigen} and Lemma \ref{lem:covdensity} are not included since they are stand-alone results, giving conditions under which  bounded eigenvalue assumptions hold. \\

%\addtocounter{section}{-1}

\section{Alternative Parametrization}
\label{sec:model2}

In the main text we focused on the model that puts smoothness assumptions on $(\pi,\mu_0,\tau)$, and gave a detailed proof of the minimax lower bound in the $\alpha \geq \beta$ case. In this section we give lower bound constructions for the $\beta \geq \alpha$ case, as well as minimax lower bounds for a different model that puts smoothness assumptions on $(\pi,\eta,\tau)$, i.e., the marginal regression $\eta(x)$ rather than the control regression $\mu_0(x)$. \\

In the main text we also gave an estimator that was minimax optimal for all $(\alpha,\beta)$ in the $(\pi,\mu_0,\tau)$ model (under some conditions on the covariate density). This estimator is also optimal in the $(\pi,\eta,\tau)$ model, as described shortly in Section \ref{subsec:model2ub}. Here we also give matching upper bounds for the $(\pi,\eta,\tau)$ model using a different, adapted but similar higher-order estimator (which itself is also  minimax optimal in the $\alpha  \geq \beta$ case in the $(\pi,\mu_0,\tau)$ model). \\

\subsection{Lower Bounds}
\label{subsec:model2lb}

The following table summarizes the constructions used to derive all minimax lower bounds stated in this paper. \\

\begin{table}[h!]
\centering
\resizebox{.75\textwidth}{!}{
\begin{tabular}{ |ll|cc| } 
 \hline
 & & \multicolumn{2}{|c|}{ Construction} \\
 Model & Regime & $p_\lambda$ & $q_\lambda$ \\ 
 \hline
$\mu_0$ is $\beta$-smooth & $\alpha \geq \beta$ & $(1/2,\mu_{0\lambda}-\tau_h/2, \tau_h)$ & $(\pi_\lambda,\mu_{0\lambda},0)$ \\ 
& $\beta \geq \alpha$ & $(\pi_\lambda,1/2-\tau_h/2,\tau_h)$ & $(\pi_\lambda,\mu_{0\lambda},0)$ \\ 
$\eta$ is $\beta$-smooth & $\alpha \geq \beta$ & $(1/2,\eta_\lambda,\tau_h)$ & $(\pi_\lambda,\eta_\lambda,0)$ \\ 
 & $\beta \geq \alpha$ & $(1/2,1/2,1-\tau_h)$ & $(\pi_\lambda,1/2,1)$ \\ 
 \hline
\end{tabular} }
\vspace{.1in}
\caption{The constructions used for our lower bound arguments. For the $\mu_0$ parametrization the density components are $(\pi,\mu_0,\tau)$ and for the $\eta$ parametrization the components are $(\pi,\eta,\tau)$. The first row is the construction detailed in the main text. The constructions in the first and third rows are equivalent, since $\eta=\pi\tau+\mu_0$, so that $\eta=\mu_{0\lambda}$ in the first row.  These first/third rows are the construction used in the first preprint of this paper.}
 \label{tab:constructions}
\end{table}

Theorem \ref{thm:lowerbound2} below states the minimax rate for the $(\pi,\eta,\tau)$ model that puts smoothness assumptions on the marginal regression function $\eta$, rather than the control regression $\mu_0$. \\

\begin{theorem} \label{thm:lowerbound2}
For $x_0 \in (0,1)$, let $\mathcal{P}$ denote the model where:
\begin{enumerate}
\item $f(x)$ is bounded above by a constant,
\item $\pi(x)$ is $\alpha$-smooth, %and $\epsilon \leq \pi(x) \leq 1-\epsilon$ for some $\epsilon>0$, 
\item $\eta(x)$ is $\beta$-smooth, and
\item $\tau(x)$ is $\gamma$-smooth.
\end{enumerate}
Let $s \equiv (\alpha+\beta)/2$. 
Then for $n$ larger than a constant depending on $(\alpha,\beta,\gamma,d)$, the minimax rate is lower bounded as
\begin{align*}
\inf_{\widehat\tau} \sup_{P \in \mathcal{P}} \E_P | \widehat\tau(x_0) - \tau_P(x_0)| \gtrsim \begin{cases}
n^{-1/\left(1 + \frac{d}{2\gamma} + \frac{d}{4 \min(\alpha,s)} \right) } & \text{ if } \min(\alpha,s) < \frac{d/4}{1+d/2\gamma}  \\
n^{-1/\left(2 + \frac{d}{\gamma}  \right) } & \text{ otherwise}  . 
\end{cases}
\end{align*}
\end{theorem}

\bigskip

 In the following subsections we give proofs of Theorems \ref{thm:lowerbound} and \ref{thm:lowerbound2} for the $\beta \geq \alpha$ regimes not detailed in the main text. Given the constructions listed in Table \ref{tab:constructions}, the arguments are all very similar to those presented in the main text for the $\alpha \geq \beta$ case. \\

\subsubsection{Setting 1: $\alpha \geq \beta$}

When $\alpha \geq \beta$, the same construction can be used for both the $\mu_0$ and $\eta$ parametrizations. This is the construction in the first and third rows of Table \ref{tab:constructions}, and detailed in the main text. Thus the proofs there settle this case. \\

\subsubsection{Setting 2: $\beta \geq \alpha$ and $\eta$ is $\beta$-smooth}

For the construction in the 4th row of Table \ref{tab:constructions}, $\eta(x)=1/2$ is constant and so certainly $\beta$-smooth, and $\pi(x)$ is $\alpha$-smooth as in the main text. Further we have by Proposition \ref{prop:dens} that
\begin{align*}
p_\lambda(z) &= f(x)  \left[ \frac{1}{4}   + (2a-1)(2y-1) \left\{ \frac{1}{4} - \frac{h^\gamma}{4} B\left( \frac{x-x_0}{2h} \right) \right\} \right] \\
q_\lambda(z) &= f(x)  \left[ \frac{1}{4} + (a-1/2) \Delta_\pi(x)  + (2a-1)(2y-1)  \left\{ \frac{1}{4} - \Delta_\pi(x)^2 \right\}  \right]
\end{align*}
with $\Delta_\pi(x)= \left(h/k^{1/d} \right)^\alpha \sum_{j=1}^k \lambda_j B\left( \frac{x-m_j}{h/k^{1/d}} \right)$, so that $\overline{p}(z) = p_\lambda(z)$ and
\begin{align*}
\overline{q}(z) &= f(x)  \left[ \frac{1}{4}  + (2a-1)(2y-1)  \left\{ \frac{1}{4} - \left(h/k^{1/d} \right)^{2\alpha} \sum_{j=1}^k B\left( \frac{x-m_j}{h/k^{1/d}} \right)^2 \right\}  \right] .
\end{align*}
Now we detail the quantities $(\delta_1,\delta_2,\delta_3)$ from Lemma \ref{lem:hellbound}, which yield a bound on the Hellinger distance between mixtures of $n$-fold products from $P_\lambda$ and $Q_\lambda$. \\

First, since $\overline{p}(z) = p_\lambda(z)$ it immediately follows that $\delta_1=0$. Similarly, taking 
$$ \frac{h^\gamma}{4} = \left(h/k^{1/d} \right)^{2\alpha} $$
ensures that $\overline{p}(z)=\overline{q}(z)$ and thus $\delta_3=0$, just as in \eqref{eq:delta3} in the main text. 
For $\delta_2$, note that  $q_\lambda(z) - p_\lambda(z) = f(x)(a-1/2) \Delta_\pi(x)$ (with the above choices of $h$ and $k$)
so that
$$ \delta_2 \leq \left( \frac{2^d \|B \|_2^2}{\epsilon/2} \right)  \left( h/k^{1/d} \right)^{2\alpha} $$
by the same arguments as for $\delta_1$ in the main text. Now by Lemma \ref{lem:hellbound} the Hellinger distance is bounded above as
$$ H^2\left(\int P_\lambda^n \ d\varpi(\lambda), \int Q_\lambda^n \ d\varpi(\lambda) \right) \ \lesssim \  \left( \frac{n^2 h^{2d}}{k} \right)  \left( h/k^{1/d} \right)^{4\alpha} $$  
so taking 
$$ \frac{h}{k^{1/d}} = \left( \frac{h^\gamma}{4} \right)^{1/2\alpha} \sim \left( \frac{1}{n^2} \right)^{\frac{1}{d+ 4\alpha+2\alpha d/\gamma}} $$
ensures the Hellinger distance is bounded above by one, and yields the promised minimax rate, since the separation in $\tau(x_0)$ under this construction is
$$ \frac{h^\gamma}{4} = \left( \frac{h}{k^{1/d}} \right)^{2\alpha} \sim n^{-1/\left(1 + \frac{d}{2\gamma} + \frac{d}{4 \alpha} \right) } .  $$

\bigskip

\subsubsection{Setting 3: $\beta \geq \alpha$ and $\mu_0$ is $\beta$-smooth}

For the construction in the 2nd row of Table \ref{tab:constructions} it is again clear that $\pi(x)$ and $\mu_0(x)$ are $\alpha$- and $\beta$-smooth, respectively. Note for the latter that, under $p_\lambda$, the control regression $\mu_0(x) = 1/2-\tau_h(x)/2$ is $\gamma$-smooth and $\gamma \geq \beta$ by definition. 
By Proposition \ref{prop:dens} we have
$$ q_\lambda(z) = f(x) \left\{ \frac{1}{4} + (a-1/2) \Delta_\pi(x) + (y-1/2)  \Delta_\mu(x)  + (2a-1)(2y-1)  \Delta_\pi(x) \Delta_\mu(x)   \right\} $$
with $\Delta_\pi(x)= \left(h/k^{1/d} \right)^\alpha \sum_{j=1}^k \lambda_j B\left( \frac{x-m_j}{h/k^{1/d}} \right)$ and $\Delta_\mu(x)= \left(h/k^{1/d} \right)^\beta \sum_{j=1}^k \lambda_j B\left( \frac{x-m_j}{h/k^{1/d}} \right)$, as in the main text, and 
\begin{align*}
p_\lambda(z) &=f(x) \bigg( \frac{1}{4} + (a-1/2) \Delta_\pi(x) - (y-1/2)  \tau_h(x)/2  \\
& \hspace{.5in} + (2a-1)(2y-1) \left[ -\Delta_\pi(x) \tau_h(x)/2 + a \tau_h(x) \left\{ \frac{1}{2} + \Delta_\pi(x) \right\}   \right] \bigg) \\
%&= f(x) \bigg( \frac{1}{4} + (a-1/2) \Delta_\pi(x) + (2a-1)(2y-1) \\
%& \hspace{.5in} \times \left[ -\Delta_\pi(x) \tau_h(x)/2 - (2a-1) \tau_h(x)/4 + a \tau_h(x) \left\{ \frac{1}{2} + \Delta_\pi(x) \right\}   \right] \\
&= f(x) \left[ \frac{1}{4} + (a-1/2) \Delta_\pi(x) + (2a-1)(2y-1) \left\{ \frac{\tau_h(x)}{4} + (a-1/2) \Delta_\pi(x) \tau_h(x) \right\} \right]
\end{align*}
where the second equality follows since $(2a-1)^2=1$. Therefore
\begin{align*}
\overline{p}(z) &= f(x) \left\{ \frac{1}{4} + (2a-1)(2y-1) \frac{h^\gamma}{4} B\left( \frac{x-x_0}{2h} \right) \right\} \\
\overline{q}(z) &= f(x) \left\{ \frac{1}{4} + (2a-1)(2y-1) \left( h/k^{1/d} \right)^{2s} \sum_{j=1}^k B\left( \frac{x-m_j}{h/k^{1/d}} \right)^2 \right\}  
\end{align*}
as in the main text. 
Now we again detail the quantities $(\delta_1,\delta_2,\delta_3)$ from Lemma \ref{lem:hellbound}. \\

Since the marginalized densities $\overline{p}(z)$ and $\overline{q}(z)$ here take the same form as those in the main text, the same arguments imply $\delta_3=0$ as long as $h^\gamma/4 =(h/k^{1/d})^{2s}$. 
Since
\begin{align*}
\{ p_\lambda(z) - \overline{p}(z) \}^2 &= \left[ f(x) \Big\{ (a-1/2) + (y-1/2)\tau_h(x) \Big\} \Delta_\pi(x) \right]^2 \leq f(x)^2 \Delta_\pi(x)^2
\end{align*}
as long as $\tau_h(x) \leq 1$, we have that
$$ \delta_1 \lesssim \left( h/k^{1/d} \right)^{2\alpha} $$
by the same arguments as for $\delta_1$ in the main text (and $\delta_2$ in the previous subsection). Similarly, with the above choices of $h$ and $k$ ensuring $\delta_3=0$, we have
\begin{align*}
\{ q_\lambda(z) - p_\lambda(z) \}^2 &= \left[ f(x) (y-1/2) \Big\{ \Delta_\mu(x) - \tau_h(x) \Delta_\pi(x) \Big\} \right]^2 \\
&=  \left[ f(x) (y-1/2) \Delta_\mu(x) \Big\{ 1 - 4 \Delta_\pi(x)^2 \Big\} \right]^2 \\
&\lesssim f(x)^2 \Delta_\mu(x)^2
\end{align*}
where the second equality follows since $\tau_h(x)/4=\Delta_\pi(x)\Delta_\mu(x)$, and the last inequality since $\Delta_\pi(x)^2$ is bounded. Therefore
$$ \delta_2 \lesssim  \left( h/k^{1/d} \right)^{2\beta}   $$
by the same arguments as for $\delta_1$ and $\delta_2$ in the main text. This implies by Lemma \ref{lem:hellbound} that the Hellinger distance is bounded above by
$$ \left( \frac{n^2 h^{2d}}{k} \right)  \left\{ \left( h/k^{1/d} \right)^{4s} + \left( h/k^{1/d} \right)^{4\beta}  \right\} \ \lesssim \ \left( \frac{n^2 h^{2d}}{k} \right)   \left( h/k^{1/d} \right)^{4s} $$  
where the inequality follows since $\beta \geq \alpha$ in this setting. Now, as in the main text, taking 
$$ \frac{h}{k^{1/d}} = \left( \frac{h^\gamma}{4} \right)^{1/2s} \sim \left( \frac{1}{n^2} \right)^{\frac{1}{d+ 4s+2s d/\gamma}} $$
ensures the Hellinger distance is bounded above by one, and yields the promised minimax rate, since the separation in $\tau(x_0)$ under this construction is
$$ \frac{h^\gamma}{4} = \left( \frac{h}{k^{1/d}} \right)^{2s} \sim n^{-1/\left(1 + \frac{d}{2\gamma} + \frac{d}{4 s} \right) } .  $$

\bigskip

\subsection{Upper Bounds}
\label{subsec:model2ub}

In the model where $\eta$ instead of $\mu_0$ is assumed to be $\beta$-smooth, the estimator described in the main text is also minimax optimal (since if $\eta$ is $\beta$-smooth, then $\mu_0=\eta - \pi \tau  $ is $\min(\alpha,\beta)$-smooth, and so in the $\alpha < \beta$ case one gets the same rates in the main text with $s$ replaced by $\min(s,\alpha)$). However one can also use the following alternative estimator, which we detail for posterity. \\

\begin{definition}[Higher-Order Local Polynomial R-Learner] \label{def:estimator2} 
Let $K_h(x)$ and $\rho_h(x)$ be defined exactly as in Definition \ref{def:estimator}. 
The proposed estimator is then defined as
\begin{equation} \label{eq:estimator2}
\widehat\tau(x_0) = \rho_h(x_0)^\T \widehat{Q}^{-1} \widehat{R} 
\end{equation}
where $\widehat{Q}$ is a $q \times q$ matrix and $\widehat{R}$ a $q$-vector now given by
\begin{align*}
\widehat{Q} &= \Pn \Big\{  \rho_h(X) K_h(X)\widehat\varphi_{a1}(Z)  \rho_h(X)^\T  \Big\} + \Un \Big\{ \rho_h(X_1) K_h(X_1) \widehat\varphi_{a2}(Z_1,Z_2) K_h(X_2) \rho_h(X_1)^\T \Big\} \\
\widehat{R} &=  \Pn \Big\{  \rho_h(X_1) K_h(X_1)  \widehat\varphi_{y1}(Z_1) \Big\} + \Un \Big\{  \rho_h(X_1) K_h(X_1) \widehat\varphi_{y2}(Z_1,Z_2) K_h(X_2) \Big\} , 
%\widehat{Q} &= \Un\left[ \rho(X_1) K_h(X_1) \Big\{ \widehat\varphi_{a1}(Z_1) + \widehat\varphi_{a2}(Z_1,Z_2) K_h(X_2) \Big\} \rho(X_1)^\T \right] \\
%\widehat{R} &=  \Un \left[ \rho(X_1) K_h(X_1) \Big\{ \widehat\varphi_{y1}(Z_1) +  \widehat\varphi_{y2}(Z_1,Z_2) K_h(X_2) \Big\} \right]
\end{align*}
respectively, with 
\begin{align*}
\widehat\varphi_{a1}(Z) &= \{A-\widehat\pi(X)\}^2 \\
\widehat\varphi_{y1}(Z) &= \{Y - \widehat\eta(X) \} \{A-\widehat\pi(X)\} \\
\widehat\varphi_{a2}(Z_1,Z_2) &= -\{A_1 - \widehat\pi(X_1)\} b_{hk}(X_1)^\T \widehat\Omega^{-1} b_{hk}(X_2) \{A_2 - \widehat\pi(X_2)\} \\
\widehat\varphi_{y2}(Z_1,Z_2) &= -\{A_1 - \widehat\pi(X_1)\} b_{hk}(X_1)^\T \widehat\Omega^{-1} b_{hk}(X_2) \{Y_2 - \widehat\eta(X_2)\}  \\
b_{hk}(x) &=  b\{1/2 + (x-x_0)/h\} \one(\|x-x_0 \| \leq h/2) \\
\widehat\Omega &= \int_{v \in [0,1]^d} b(v) b(v)^\T \ d\widehat{F}(x_0 + h(v-1/2)) 
\end{align*}
for $b: \R^d \mapsto \R^k$ a basis of dimension $k$. 
The nuisance estimators $( \widehat{F}, \widehat\pi,\widehat\eta)$ are constructed from a separate training sample $D^n$, independent of that on which $\Un$ operates. \\
\end{definition}

Like the estimator \eqref{eq:estimator} in the main text, the estimator \eqref{eq:estimator2}  can also be viewed as a modified second-order estimator of the projection parameter. However, instead of estimating the moment condition \eqref{eq:moment} under a linear effect-type assumption, the estimator \eqref{eq:estimator2} instead estimates the numerator and denominator terms $R$ and $Q$ separately. For example, the first term in $\widehat{R}$, i.e., 
$\Pn \{ \rho_h(X) K_h(X)  \widehat\varphi_{y1}(Z) \} , $
is the usual first-order influence function-based estimator of $R$. \citet{kennedy2023towards} analyzed an undersmoothed version of this  estimator (where the nuisance estimates $\widehat\pi$ and $\widehat\mu$ themselves are undersmoothed linear smoothers), calling it the local polynomial R-learner. 
The second term 
$$ \Un  \Big\{  \rho_h(X_1) K_h(X_1) \widehat\varphi_{y2}(Z_1,Z_2) K_h(X_2) \Big\} $$
is similar to the second-order U-statistic correction that would be added using the higher-order influence function methodology developed by \citet{robins2008higher, robins2009quadratic, robins2017minimax}. However,  this term differs in its localization just as described for the estimator in the main text. \\

\begin{propositionapp} \label{prop:qrerror2}
The estimator \eqref{eq:estimator2} satisfies
$$ \left| \widehat\tau(x_0) - \tau_h(x_0) \right| \ \leq \ \|\rho(1/2) \| \left(|  \| Q^{-1} \| +   \| \widehat{Q}^{-1} - {Q}^{-1} \| \right) \left( \left\| \widehat{R} - R \right\|  + \left\|  Q - \widehat{Q} \right\|_2 \| Q^{-1} R \| \right) , $$
and further if 
 $\| Q^{-1} \|$, $ \| \widehat{Q}^{-1} - {Q}^{-1} \|$, and $\|Q^{-1} R\| $ are all bounded above, 
then
\begin{align*}
\E \left| \widehat\tau(x_0) - \tau_h(x_0) \right| \ &\lesssim \  \max_j  \sqrt{ \E\left\{ \E(\widehat{R}_j - R_j \mid D^n)^2 + \var(\widehat{R}_j \mid D^n) \right\} }   \\
& \hspace{.5in} + \max_{j,\ell} \sqrt{ \E\left\{ \E(\widehat{Q}_{j\ell} - Q_{j\ell} \mid D^n)^2 + \var(\widehat{Q}_{j\ell} \mid D^n) \right\} } .
\end{align*}
for $D^n$  a separate independent training sample on which $( \widehat{F}, \widehat\pi,\widehat\eta)$ are estimated. \\
\end{propositionapp}

\bigskip

\begin{proof}
We have
\begin{align*}
| \widehat\tau(x_0) - \tau_h(x_0) | %&= \left| \rho(x_0)^\T \left( \widehat{Q}^{-1} \widehat{R} - Q^{-1} R \right) \right| \\
&= \left| \rho_h(x_0)^\T \widehat{Q}^{-1} \left\{  \left( \widehat{R} - R \right) +    \left( Q - \widehat{Q} \right) Q^{-1} R \right\} \right| \\
&\leq \|\rho(1/2) \| \left\| \widehat{Q}^{-1}  \right\| \left( \left\| \widehat{R} - R \right\|  + \left\|  Q - \widehat{Q} \right\| \left\| Q^{-1} R \right\| \right) \\
& \leq \|\rho(1/2) \| \left( \left\| Q^{-1} \right\|+ \left\| \widehat{Q}^{-1} - Q^{-1} \right\| \right) \left( \left\| \widehat{R} - R \right\|  + \left\|  Q - \widehat{Q} \right\|_2 \left\| Q^{-1} R \right\| \right)
\end{align*}
by the sub-multiplicative and triangle inequalities of the operator norm, along with the fact that $\| A \| \leq \| A \|_2$. Together with the bounds on $\|Q^{-1} \|$, $\| \widehat{Q}^{-1} - Q^{-1} \|$, and $\| Q^{-1} R \|$, this yields the first inequality. For the second inequality, first note $\| \rho(x)\| \leq Cq$, as described in the proof of Proposition \ref{prop:projbias}.  The second inequality now  follows since
\begin{align*}
\E\| \widehat{R} - R  \|  &\leq \sqrt{ \E\| \widehat{R} - R  \|^2 } %= \sqrt{ \sum_j \E\Big\{ (\widehat{R}_j - R_j)^2 \Big\}  } \\
= \sqrt{ \sum_j \E\left[ \E\Big\{ (\widehat{R}_j - R_j)^2 \mid D^n \Big\} \right] } \\
&= \sqrt{ \sum_j \E\left\{ \bias(\widehat{R}_j \mid D^n)^2 + \var(\widehat{R}_j \mid D^n) \right\} } \\
&\leq \sqrt{ {{ d + \lfloor \gamma \rfloor } \choose { \lfloor\gamma \rfloor }} } \max_j \sqrt{ \E\left\{ \bias(\widehat{R}_j \mid D^n)^2 + \var(\widehat{R}_j \mid D^n) \right\} } 
\end{align*}
using Jensen's inequality and iterated expectation. The last line follows since the length of $R$ is ${{ d + \lfloor \gamma \rfloor } \choose { \lfloor\gamma \rfloor }}$. The logic is the same for $\E\| \widehat{Q} - Q \|_2 = \E \sqrt{ \sum_{j, \ell}  (\widehat{Q}_{j\ell} - Q_{j\ell})^2}$. 
\end{proof}

\bigskip

\begin{propositionapp} \label{prop:bias_alt}
%Let $\| f \|_{K_h}^2 = \int f(x)^2 K_h(x) \ dF(x)$ denote the (squared) kernel-weighted $L_2(\Pb)$ norm. 
Assume:
\begin{enumerate}
\item $\lambda_{\max}(\Omega)$ is bounded above, 
%\item $\|dF^*/dF\|_\infty$ is bounded above, where $F^*(v) = F(x_0 + h(v-1/2))$, 
\item the basis $b$ satisfies approximating condition \eqref{eq:holder}, 
\item $\widehat\pi(x) - \pi(x)$ is $\alpha$-smooth, 
\item $\widehat\eta(x) - \eta(x)$ is $\beta$-smooth.
\end{enumerate}
Then
\begin{align*}
| \E(\widehat{R}_j - R_j \mid D^n) | & \ \lesssim \ \left(h/k^{1/d} \right)^{2s} + \| \widehat\pi - \pi \|_{F^*} \| \widehat\eta - \eta \|_{F^*} \| \widehat\Omega^{-1}-\Omega^{-1} \|  \\
| \E(\widehat{Q}_{j\ell} - Q_{j\ell} \mid D^n) | & \ \lesssim \ \left( h/k^{1/d} \right)^{2\alpha} + \| \widehat\pi - \pi \|_{F^*}^2 \| \widehat\Omega^{-1}-\Omega^{-1} \| . 
\end{align*}
\end{propositionapp}

\bigskip

\begin{proof}
We only prove the result for $\widehat{R}_j$, since the logic is the same for $\widehat{Q}_{j\ell}$. By iterated expectation, the conditional mean of the first-order term is
\begin{align} 
\E\{ \rho_h(X_1)  K_h(X_1) & \widehat\varphi_{y1}(Z) \mid D^n \} = R + \int  \rho_h(x) K_h(x)  \{ \pi(x) - \widehat\pi(x) \}\{ \eta(x) - \widehat\eta(x) \}  \ dF(x) \nonumber \\
& = R + \int \rho(v) \{ \pi^*(v) - \widehat\pi^*(v) \}\{ \eta^*(v) - \widehat\eta^*(v) \}  \ dF^*(v) \label{eq:firstbias2}
\end{align}
where we use the change of variable $v = \frac{1}{2} + \frac{x-x_0}{h}$ and again define for any function $g: \R^d \mapsto \R$ its corresponding stretched version as $g^*(v) = g(x_0 + h(v-1/2))$. To ease notation it is left implicit that any integral over $v$ is only over $\{v : \|v-1/2\| \leq 1/2 \}$. 
Similarly, (minus) the conditional mean of the second-order influence function term is
\begin{align}
&-\E\{   \rho_h(X_1) K_h(X_1) \widehat\varphi_{y2}(Z_1,Z_2)K_h(X_2) \mid D^n \} \nonumber \\
&=  \iint  \rho_h(x_1) K_h(x_1) \{\pi(x_1) - \widehat\pi(x_1)\} b_{hk}(x_1)^\T \widehat\Omega^{-1} b_{hk}(x_2) \{ \eta(x_2) - \widehat\eta(x_2)\} K_h(x_2) \ dF(x_2) \ dF(x_1) \nonumber \\
&=  \iint \rho(v_1)  \{\pi^*(v_1) - \widehat\pi^*(v_1)\} b(v_1)^\T \widehat\Omega^{-1} b(v_2) \{ \eta^*(v_2) - \widehat\eta^*(v_2)\}  \ dF^*(v_2) \ dF^*(v_1) \nonumber \\
&=  \int \rho(v_1)\{\pi^*(v_1) - \widehat\pi^*(v_1)\} \widehat\Pi_{b}(\eta^*-\widehat\eta^*)(v_1) \ dF^*(v_1) \label{eq:secondbias2}
\end{align}
where we define
\begin{align*}
\Pi_{b}g^*(u) &= b(u)^\T \Omega^{-1} \int {b}(v) g^*(v) \ dF^*(v) 
\end{align*}
as the $F^*$-weighted linear projection of $g^*$ on the basis $b$, and $\widehat\Pi_{b} g^*(u)$ as the estimated version, which simply replaces $\Omega$ with $\widehat\Omega$. 
Therefore adding the first- and second-order expected values in \eqref{eq:firstbias2} and \eqref{eq:secondbias2}, the overall bias relative to $R$ is
\begin{align}
\int  \rho(v) & \{ \pi^*(v) - \widehat\pi^*(v)\} (I-\widehat\Pi_{b}) (\eta^* - \widehat\eta^*)(v) \ dF^*(v) \nonumber \\
%&= \int \rho(v) \{ \pi^*(v) - \widehat\pi^*(v)\} (I-\Pi_{b}) (\eta^* - \widehat\eta^*)(v) \ dF^*(v) \nonumber \\
%& \hspace{.6in} + \int \rho(v) \{ \pi^*(v) - \widehat\pi^*(v)\} (\Pi_{b}-\widehat\Pi_{b}) (\eta^* - \widehat\eta^*)(v) \ dF^*(v)  \nonumber \\
&= \int (I-\Pi_b) \{ \rho (\pi^* - \widehat\pi^*)\}(v) (I-\Pi_{b}) (\eta^* - \widehat\eta^*)(v) \ dF^*(v) \label{eq:mainbias2}  \\
& \hspace{.6in} + \int \rho(v) \{ \pi^*(v) - \widehat\pi^*(v)\} (\Pi_{b}-\widehat\Pi_{b}) (\eta^* - \widehat\eta^*)(v) \ dF^*(v)  \label{eq:covbias2} 
\end{align}
where the last line follows from orthogonality of a projection with its residuals (Lemma \ref{lem:projections}(i)). \\

Now we analyze the bias terms \eqref{eq:mainbias2} and \eqref{eq:covbias2} separately; the first is the main bias term, which would arise even if the covariate density were known, and the second is the contribution coming from having to estimate the covariate density. \\

Recall from \eqref{eq:smallbias} in the main text that if $g$ is $\alpha$-smooth, then 
$\| (I - \Pi_b) g^*/h^\alpha \|_{F^*} \lesssim k^{-\alpha/d}$
by the \Holder{} approximation properties \eqref{eq:holder} of the basis $b$, and so it follows that
\begin{equation*}
\| (I - \Pi_b) g^* \|_{F^*} \ \lesssim \ h^\alpha k^{-\alpha/d} = (k/h^d)^{-\alpha/d} 
\end{equation*}
for any $\alpha$-smooth function $g$. Therefore now consider the bias term \eqref{eq:mainbias2}. This term satisfies
\begin{align*}
\int \Big[ &(I- \Pi_{b})\{\rho(\pi^*-\widehat\pi^*)\}(v) \Big] \Big\{ (I-\Pi_{b}) (\eta^* - \widehat\eta^*)(v) \Big\} \ dF^*(v) \\
&\leq \| (I-\Pi_{b}) \{\rho(\pi^*-\widehat\pi^*)\}  \|_{F^*} \| (I-\Pi_{b})(\eta^*-\widehat\eta^*) \|_{F^*}  \\
&  \lesssim  \  (k/h^d)^{-2s/d}
\end{align*}
where the second line follows by Cauchy-Schwarz, and the third by \eqref{eq:smallbias}, 
since $(\pi-\widehat\pi)$ and  $(\eta-\widehat\eta)$ are assumed $\alpha$- and $\beta$-smooth, respectively (note $\rho(v)$ is a polynomial, so the smoothness of $\rho(\pi^*-\widehat\pi^*)$ is the same as $(\pi^*-\widehat\pi^*)$). \\

Now for the term in  \eqref{eq:covbias2}, let $\theta_{b,g} = \Omega^{-1} \int b g \ dF^*$ denote the coefficients of the projection $\Pi_{b}g$, and note for any functions $g_1,g_2$ we have
\begin{align*}
\int g_1(\Pi_{b} - \widehat\Pi_{b})(g_2) \ dF^* %\\ &=  \left\{ \int {b}_h(x) K_h(x) g(x)   \ dF(x) \right\}^\T   (\Omega^{-1} - \widehat\Omega^{-1})  \left\{   \int {b}_h(t) K_h(t) f(t)  \ d\Pb(t) \right\} \\
&= \left( \Omega^{1/2} \theta_{b, g_1} \right)^\T \Omega^{1/2} (\Omega^{-1} - \widehat\Omega^{-1}) \Omega^{1/2}  \left( \Omega^{1/2} \theta_{b, g_2} \right) \\
&\leq \| g_1 \|_{F^*}  \| \Omega^{1/2} (\Omega^{-1} - \widehat\Omega^{-1}) \Omega^{1/2}  \| \| g_2 \|_{F^*} \\
%&= \| f \|_{K_h} \| g \|_{K_h}  \| \Omega^{1/2} \Omega^{-1}  ( \widehat\Omega - \Omega)  \widehat\Omega^{-1} \Omega^{1/2}  \| \\
&\leq  \| g_1 \|_{F^*} \| g_2 \|_{F^*}  \|  \Omega \| \| \widehat\Omega^{-1} - \Omega^{-1} \|  
\end{align*}
where the first equality follows by definition, the second line since the $L_2$ norm of the coefficients of a (weighted) projection is no more than the weighted $L_2(\Pb)$ norm of the function itself (Lemma \ref{lem:projections}(iii)), and the last by the sub-multiplicative property of the operator norm, along with the fact that $\| \Omega^{1/2} \|^2 = \| \Omega \|$. 
\end{proof}

\bigskip

We omit details on variance calculations because the logic is exactly the same as for the estimator \eqref{eq:estimator} studied in Section \ref{sec:variance} of the main text. \\

Combining the approximation bias in Proposition \ref{prop:projbias} with the decomposition in Proposition \ref{prop:qrerror2}, and the relevant bias and variance bounds together with sufficient conditions on covariate density estimation, we obtain that the estimator \eqref{eq:estimator2} has error satisfying
\begin{equation} \label{eq:overallrate2}
 \E_P | \widehat\tau(x_0) - \tau_P(x_0) |  \ \lesssim \ h^\gamma + \left( h/k^{1/d} \right)^{2\min(s,\alpha)} + \sqrt{\frac{1}{nh^d} \left( 1 + \frac{k}{nh^d} \right) } 
 \end{equation}
Balancing terms as in Section \ref{sec:overallrate} yields the minimax rate from Theorem \ref{thm:lowerbound2}. \\

\section{Additional Results}
\label{sec:propproofs}

\subsection{Density Bounds in Proposition \ref{prop:hellbound}}
\label{app:prophellbound}

Here we show that the relevant densities and density ratios from Proposition \ref{prop:hellbound} are appropriately bounded. \\

\begin{proof} 
When $h \leq 1/4$ then it follows that on $\mathcal{S}_{hk}$ we have
\begin{equation}
1 \leq f(x) = \left\{ 1  - \left( \frac{4^d - 1}{2^d} \right) h^d \right\}^{-1} \leq 2 . \label{eq:fbd}
\end{equation}
Further, since $B(x) \leq \one(x \in [-1,1]^d)$, $a,y \in \{0,1\}$, and $\lambda \in \{-1,1\}$, we have on $\mathcal{S}_{hk}$ that 
\begin{equation*}
\frac{1}{4} - \frac{ 1 }{2} \left( h/k^{1/d} \right)^\beta  - \frac{h^\gamma}{4}   \  \leq \frac{p_\lambda(z)}{f(x)}  \leq  \  \frac{1}{4} + \frac{ 1}{2}  \left( h/k^{1/d} \right)^\beta  + \frac{h^\gamma}{4}   , 
\end{equation*}
regardless of the values of $h,k \geq 0$. Therefore when $h^\gamma + 2(h/k^{1/d})^\beta \leq 1- 4\epsilon$, the above bound implies
\begin{equation}
\frac{p_\lambda(z)}{f(x)}  \geq  \frac{1}{4} - \frac{ 1 }{2}  \left( h/k^{1/d} \right)^\beta - \frac{h^\gamma}{4}  \geq \epsilon . \label{eq:plamfbd}
\end{equation}
Similarly, when $h^\gamma + 2(h/k^{1/d})^\beta \leq 1- 4\epsilon$ (which implies $h^\gamma \leq 1$) we also have
\begin{equation*}
\frac{\overline{p}(z)}{p_\lambda(z)} \leq \frac{\frac{1}{4} + \frac{h^\gamma}{4} }{\frac{1}{4}  - \frac{ 1 }{2}  \left( h/k^{1/d} \right)^\beta - \frac{h^\gamma}{4} } \leq \frac{1/2}{\epsilon}  .
\end{equation*}
Note that, although \citet{robins2009semiparametric} assume $p_\lambda(z)$ is uniformly lower bounded away from zero in their version of Lemma \ref{lem:hellbound}, they only use a bound on $\overline{p}/p_\lambda$ to ensure their quantity $c$ is bounded (see page 1319). Therefore this condition also holds in our case. 
\end{proof}

\bigskip

\subsection{Density Expressions} 

\begin{propositionapp} \label{prop:dens}
For binary $A$ and $Y$, the density of $Z=(X,A,Y)$ can be written as
\begin{align*}
f(x) & \bigg( \frac{1}{4} + (a-1/2) \Delta_\pi(x) + (y-1/2)  \Delta_\mu(x)  \\
& \hspace{.4in}  + (2a-1)(2y-1) \left[ \Delta_\pi(x) \Delta_\mu(x) + a \tau(x) \left\{ \frac{1}{2} + \Delta_\pi(x) \right\}   \right] \bigg) 
\end{align*}
for $\Delta_\pi(x) = \pi(x) - 1/2$ and $\Delta_\mu(x) = \mu_0(x)-1/2$, or equivalently as
\begin{align*}
f(x) & \bigg( \frac{1}{4} + (a-1/2) \Delta_\pi(x) + (y-1/2)  \Delta_\eta(x)  \\
&\hspace{.4in} + (2a-1)(2y-1) \left[ \Delta_\pi(x) \Delta_\eta(x) + \tau(x) \left\{ \frac{1}{4} - \Delta_\pi(x)^2 \right\}   \right] \bigg) 
\end{align*}
for $\Delta_\eta(x) = \eta(x)-1/2$. 
\end{propositionapp}

\bigskip

\begin{propositionapp} \label{prop:dens2}
The densities under $P_\lambda$ and $Q_\lambda$ from Definition \ref{def:construction} are given by
\begin{align*}
p_\lambda(z) &= f(x) \Bigg[ \frac{1}{4} + (y-1/2) \left( h/k^{1/d} \right)^{\beta} \sum_{j=1}^k \lambda_j B \left( \frac{x-m_j}{h/k^{1/d}}  \right) + (2a-1)(2y-1) \frac{h^\gamma}{4} B\left( \frac{x-x_0}{h} \right) \Bigg] \\
q_\lambda(z) &= f(x) \Bigg[ \frac{1}{4} + \left\{ (a-1/2)  \left( h/k^{1/d} \right)^\alpha + (y-1/2)    \left( h/k^{1/d} \right)^\beta \right\} \sum_{j=1}^k \lambda_j B \left( \frac{x-m_j}{h/k^{1/d}}  \right)  \\
& \hspace{.7in}  + (2a-1)(2y-1)  \left( h/k^{1/d} \right)^{2s} \sum_{j=1}^k B \left( \frac{x-m_j}{h/k^{1/d}}  \right) ^2  \Bigg] 
\end{align*}
where $s \equiv (\alpha+\beta)/2$. 
\end{propositionapp} 
 
 \bigskip
 
We note that the densities are both equal to $1/4$ for all $x \notin \mathcal{C}_{2h}(x_0)$ away from $x_0$, since $B \left( \frac{x-m_j}{h/k^{1/d}}  \right)=0$ for $x \notin \mathcal{C}_{h/k^{1/d}}(m_j)$ and $\mathcal{C}_{h/k^{1/d}}(m_j) \subseteq \mathcal{C}_h(x_0) \subseteq \mathcal{C}_{2h}(x_0)$, and since $B\left( \frac{x-x_0}{h} \right)=0$ for $x \notin \mathcal{C}_{2h}(x_0)$. \\

\subsection{Proof of Proposition \ref{prop:projbias}}
\label{app:projbias}

\begin{proof}
To ease notation we prove the result in the $d=1$ case but the logic is the same when $d>1$. First note that the local polynomial projection operator $Lg(x') \equiv \int g(x) w_h(x;x') \ dx$ for 
$$ w_h(x;x') \equiv \rho_h(x')^\T Q^{-1} \rho_h(x) K_h(x) \pi(x) \{1-\pi(x)\} f(x) $$
reproduces polynomials, in the sense that, for any polynomial of the form $g(x) = a^\T \rho_h(x)$, $a \in \R^q$, we have
\begin{align*}
Lg(x') &= \int \Big\{ a^\T \rho_h(x) \Big\} w_{h}(x;x') \ dx  \\
&=  \rho_h(x')^\T Q^{-1} \int \rho_h(x)  K_{h}(x) \pi(x) \{1-\pi(x)\} \rho_h(x)^\T f(x) \ dx \ a  \\
&= \rho_h(x')^\T Q^{-1} Q a = a^\T \rho_h(x') = g(x') .
\end{align*}
Then for $\|x'-x_0\| \leq C_2h$ we have that $\tau_h(x')-\tau(x')=\int w_h(x;x') \Big\{ \tau(x)  - \tau(x') \Big\} dx$ equals
\begin{align}
 \int w_h(x;x') &\left\{ \sum_{j=1}^{\lfloor \gamma \rfloor -1} \frac{D^j \tau(x')}{j!} (x-x')^j + \frac{D^{\lfloor \gamma \rfloor} \tau(x^*)}{\lfloor \gamma \rfloor!} (x-x')^{\lfloor \gamma \rfloor} \right\} dx \nonumber \\
&= 0 + \int w_h(x;x')  \left\{ \frac{D^{\lfloor \gamma \rfloor}\tau(x' + \epsilon(x-x'))-D^{\lfloor \gamma \rfloor}\tau(x')}{{\lfloor \gamma \rfloor}!} \right\} (x-x')^{\lfloor \gamma \rfloor} \ dx \nonumber \\
&\leq \int  \left| w_h(x;x') \right| \frac{C_1  \|x-x'\|^{\gamma-{\lfloor \gamma \rfloor}}}{{\lfloor \gamma \rfloor}!} \|x-x'\|^{{\lfloor \gamma \rfloor} } \ dx \nonumber \\
&\leq \frac{C_1}{{\lfloor \gamma \rfloor}!} \int  \left| w_h(x;x') \right| \left( \| x- x_0 \| + C_2 h \right)^\gamma dx \nonumber \\
&\leq \frac{C_1h^\gamma}{{\lfloor \gamma \rfloor}!} \left\{ h^d \int   \left| w_h(x_0+hu; x')   \right|  (\|u\|+C_2)^{\gamma}\ du \right\} \label{eq:projbound}
\end{align}
where the first line follows by Taylor expansion (with $x^* = x' +\epsilon(x-x')$ for some $\epsilon \in [0,1]$), the second by the polynomial reproducing property, the third since $\tau$ is $\gamma$-smooth, the fourth by the triangle inequality with $\|x_0 - x' \| \leq C_2 h$, and the last by a change of variable $u=(x-x_0)/h$ (so that $dx = h^d \ du$). \\

Now the result follows since we show the term on the right in \eqref{eq:projbound} (in brackets) is bounded under the stated assumptions. Specifically it equals
\begin{align*}
&  \int (\|u\|+C_2)^\gamma \left| \rho_h(x')^\T Q^{-1} \rho_h(x_0+hu) \right|  \one(\|u\| \leq 1/2) \\
& \hspace{.5in} \times \pi(x_0+hu) \{1-\pi(x_0+hu)\} f(x_0+hu) \ du \\
%&\leq \frac{1/4}{2^\gamma} \| \rho(x_0) \| \left\| Q^{-1} \right\|  \int_{-1/2}^{1/2} \| \rho(x_0+hu) \| f(x_0 + hu) \ du \\
&\leq (1/4)(1/2+C_2)^\gamma \|\rho(1/2) \| \left\| Q^{-1} \right\|   \int_{-1/2}^{1/2} \| \rho(1/2+u) \| f(x_0+hu) \ du \\
&=  (C_3q^2/4)(1/2+C_2)^\gamma \left\| Q^{-1} \right\| h^{-d}  \int \one(\|x-x_0\| \leq h/2) \ dF(x) \lesssim 1
\end{align*}
where the second line follows from the submultiplicative property of the operator norm and since $\pi(1-\pi) \leq 1/4$, and the last  from Assumptions 2 and 3 and since
\begin{align*}
\| \rho(x) \|^2 & \leq  {{ d + \lfloor \gamma \rfloor } \choose { \lfloor\gamma \rfloor }} (2\lfloor \gamma \rfloor  + 1) \leq C_3q^2
\end{align*}
for all $x$ (\citet{belloni2015some}, Example 3.1), since each Legendre term satisfies $| \rho_m(x_j)| \leq \sqrt{2m+1} \leq \sqrt{2\lfloor \gamma \rfloor + 1}$ for $m \leq \lfloor \gamma \rfloor$, and the length of $\rho$ is $q={{ d + \lfloor \gamma \rfloor } \choose { \lfloor\gamma \rfloor }} $, i.e., the maximum number of monomials in a polynomial in $d$ variables with degree up to $\lfloor \gamma \rfloor$. 
\end{proof}

\bigskip

\subsection{Proof of Proposition \ref{prop:qrerror}}

\begin{proof}
For the first inequality, we have
\begin{align*}
| \widehat\tau(x_0) - \tau_h(x_0) | %&= \left| \rho(x_0)^\T \left( \widehat{Q}^{-1} \widehat{R} - Q^{-1} R \right) \right| \\
&= \left| \rho_h(x_0)^\T \widehat{Q}^{-1} \left\{  \left( \widehat{R} - R \right) +    \left( Q - \widehat{Q} \right) Q^{-1} R \right\} \right| \\
&\leq \|\rho(1/2) \| \left( \left\| Q^{-1} \right\|+ \left\| \widehat{Q}^{-1} - Q^{-1} \right\| \right)\left( \left\| ( \widehat{R} - R) - (  \widehat{Q} - Q ) ( Q^{-1} R ) \right\| \right)
\end{align*}
by the sub-multiplicative and triangle inequalities of the operator norm. For the second inequality, first note $\| \rho(x)\| \leq Cq$, as described in the proof of Proposition \ref{prop:projbias}.  The result now  follows  from the bounds on $\| Q^{-1} \|$ and $ \| \widehat{Q}^{-1} - {Q}^{-1} \|$ and since
\begin{align*}
\E\| \widehat{S} - S  \|  &\leq \sqrt{ \E\| \widehat{S} - S  \|^2 } %= \sqrt{ \sum_j \E\Big\{ (\widehat{R}_j - R_j)^2 \Big\}  } \\
= \sqrt{ \sum_j \E\left[ \E\Big\{ (\widehat{S}_j - S_j)^2 \mid D^n \Big\} \right] } \\
&= \sqrt{ \sum_j \E\left\{ \bias(\widehat{S}_j \mid D^n)^2 + \var(\widehat{S}_j \mid D^n) \right\} } \\
&\leq \sqrt{ {{ d + \lfloor \gamma \rfloor } \choose { \lfloor\gamma \rfloor }} } \max_j \sqrt{ \E\left\{ \bias(\widehat{S}_j \mid D^n)^2 + \var(\widehat{S}_j \mid D^n) \right\} } 
\end{align*}
by Jensen's inequality and iterated expectation. The last line follows since $|S|={{ d + \lfloor \gamma \rfloor } \choose { \lfloor\gamma \rfloor }}$. 
\end{proof}

\bigskip

\subsection{Eigenvalues of $Q$ and $\Omega$}

\begin{propositionapp} \label{prop:qeigen}
If (i) $\epsilon \leq \pi(x) \leq 1-\epsilon$ and (ii) the density $dF(x)$ is bounded above and below away from zero on $\{x: \|x-x_0 \| \leq h/2\}$, then the eigenvalues of $Q$ are bounded above and below away from zero. 
%\begin{enumerate}
%\item $\epsilon \leq \pi(x) \leq 1-\epsilon$  and
%\item the density $dF(x)$ is bounded above and below away from zero on $\{x: \|x-x_0 \| \leq h/2\}$.
%\end{enumerate}
\end{propositionapp}

\bigskip

\begin{proof}
Define the stretched function $g^*(v) = g(x_0 + h(v-1/2))$  for any $g: \R^d \mapsto \R$. This maps values of $g$ in the small cube $[x_0-h/2,x_0+h/2]^d$ around $x_0$ to the whole space $[0,1]^d$. Then note that, with the change of variables $v=\frac{1}{2} + \frac{x-x_0}{h}$, 
\begin{align*}
Q &=  \int \rho_h(x) K_h(x) \pi(x) \{ 1-\pi(x)\} \rho_h(x)^\T \ dF(x)  \\
&= \int_{\|v -1/2\| \leq 1/2} \rho(v)  \rho(v)^\T  \pi^*(v) \{ 1-\pi^*(v)\}\ dF^*(v) .
\end{align*}
Next note that $\epsilon(1-\epsilon) \leq \pi(1-\pi) \leq 1/4$, so  the eigenvalues of $Q$ will be bounded if those of the matrix $$  \int  \rho(v) \rho(v)^\T \ dF(x_0 + h(v-1/2))  $$
are.  But $  \int  \rho(x) \rho(x)^\T \ dx = I$ by orthonormality of the Legendre polynomials on $[0,1]^d$, and  the local boundedness of $dF$ ensures $dF^*/d\mu$ is bounded above and below away from zero, for $\mu$ the uniform measure. Therefore Proposition 2.1 of  \citet{belloni2015some} yields the result. 
\end{proof}

\bigskip

\begin{propositionapp} \label{prop:oeigen}
If (i) the basis $b(v)$ is orthonormal with respect to the uniform measure, and  (ii) the density $dF(x)$ is bounded above and below away from zero on $\{x: \|x-x_0 \| \leq h/2\}$, then the eigenvalues of $\Omega$ are bounded above and below away from zero. 
\end{propositionapp}

\bigskip

\begin{proof}
The proof is similar to that of Proposition \ref{prop:qeigen}. First note that
\begin{align*}
\Omega = \int b_{hk}(x) K_h(x) b_{hk}(x)^\T dF(x) = \int b(v) b(v)^\T \ dF^*(v)
\end{align*}
by the change of variables $v=\frac{1}{2} + \frac{x-x_0}{h}$, and where $dF^*(v) = dF(x_0 + h(v-1/2))$ as before. Further $  \int b(v) b(v)^\T \ dv = I$ by the assumed orthonormality, and the local boundedness of $dF$ ensures $dF^*/d\mu$ is bounded above and below away from zero, for $\mu$ the uniform measure. Therefore Proposition 2.1 of  \citet{belloni2015some} yields the result. 
\end{proof}

\bigskip

\subsection{Proof of Proposition \ref{prop:variance}}
\label{app:propvariance}

\begin{proof}
Given the training data $D^n$,  $\widehat{S}_j=[\widehat{R}-\widehat{Q} \theta]_j$ is a second-order U-statistic with kernel
\begin{align*}
\xi(Z_1,Z_2) =  \rho_h(X_1) & K_h(X_1) \Big[ \Big\{ \widehat\varphi_{y1}(Z_1) -  \widehat\varphi_{a1}(Z_1)  \rho_h(X_1)^\T \theta \Big\} \\
& +  \Big\{ \widehat\varphi_{y2}(Z_1,Z_2)- \widehat\varphi_{a2}(Z_1,Z_2) \rho_h(X_1)^\T \theta \Big\} K_h(X_2) \Big]
\end{align*}
Thus its conditional variance is given by the usual variance of a second-order U-statistic (e.g., Lemma 6 of \citet{robins2009quadratic}), which is
\begin{align}
\var\left\{ \Un(\xi) \mid D^n \right\} %&= {n \choose k}^{-1} \sum_{c=1}^k {k \choose c} {n-k \choose k-c}  \\
&= \left( \frac{4(n-2)}{n(n-1)} \right) \var\{ \xi_1(Z_1)\mid D^n \} + \left( \frac{2}{n(n-1)} \right) \var\{ \xi(Z_1,Z_2) \mid D^n \}  \nonumber \\
&\leq  \frac{4 \E\{ \xi_1(Z_1)^2 \mid D^n \}}{n}  + \frac{4\E\{ \xi(Z_1,Z_2)^2 \mid D^n \} }{n^2} \label{eq:uvar}
\end{align}
for $\xi_1(z_1) = \int \xi(z_1,z_2) \ d\Pb(z_2)$, if $n \geq 2$. In our case $\xi_1$ equals
\begin{align*}
\rho_h(X_1) K_h(X_1) &\Big[ \widehat\varphi_{y1}(Z_1)  -  \widehat\varphi_{a1}(Z_1)  \rho_h(X_1)^\T \theta \\
&\hspace{.4in} + \{A_1 - \widehat\pi(X_1)\} \Big\{ \widehat\Pi_{b}(\mu_0^*-\widehat\mu_0^*)(X_1) - \widehat\Pi_b \widehat\pi^*(X_1) \rho_h(X_1)^\T \theta \Big\} \Big] .
\end{align*}
Therefore for the first term in \eqref{eq:uvar} we have
\begin{align*}
\int \xi_1^2 \ d\Pb &\leq 2 \bigg( \int \rho_h^2 K_h^2 (\widehat\varphi_{y1} - \widehat\varphi_{a1} \rho_h^\T\theta)^2 \ d\Pb \\
& \hspace{.4in} + \int \rho_h^2 K_h^2 \left\{ \pi(1-\pi) + (\pi -\widehat\pi)^2 \right\} \left\{ \widehat\Pi_{b}(\mu_0^*-\widehat\mu_0^*) + \widehat\Pi_b \pi^* \rho_h^\T \theta \right\}^2  dF \bigg) \\
&\lesssim \frac{1}{h^d} \left( \int \rho(v)^2 \ dF^*(v) + \int \rho(v)^2  \left[ \left\{ \widehat\Pi_{b}(\mu_0^*-\widehat\mu_0^*)(v) \right\}^2 + \left\{ \widehat\Pi_b \pi^*(v) \rho(v)^\T \theta \right\}^2 \right] \ dF^*(v) \right) \\
&\lesssim  \frac{1}{h^d} \left( 1 +  \| \widehat\mu_0 - \mu_0 \|_{F^*}^2  \right) 
\end{align*}
where the second inequality follows by the change of variables $v = \frac{1}{2} + \frac{x-x_0}{h}$, and since all of $(\widehat\varphi_{y1},\widehat\varphi_{a1},\pi,\widehat\pi)$ are uniformly bounded, $0 \leq K_h(x) \lesssim h^{-d}$, and $\sup_v \| \rho(v) \| \lesssim q$ (\citet{belloni2015some}, Example 3.1); and the third follows from Lemma \ref{lem:projections}(ii).  
Now, consider the second term in \eqref{eq:uvar}; we only detail the  term in $\xi$ involving $\widehat\varphi_{y2}$ since the logic is the same for the $\widehat\varphi_{a2}$ term, and the terms involving $(\widehat\varphi_{y1},\widehat\varphi_{a1})$ are bounded by standard arguments. Letting $\overline{b}_j(x) \equiv \Omega^{-1/2} b_j(x)$  and $M \equiv \Omega^{1/2} \widehat\Omega^{-1} \Omega^{1/2} $, we have $\E\{ \xi(Z_1,Z_2)^2 \mid D^n \}$ upper bounded by a multiple of 
\begin{align*}
%\E\{ \xi(Z_1,Z_2)^2 \mid D^n \}  %&= \int \int \Big\{ \rho(x_1) K_h(x_1) \widehat\varphi_{y2}(z_1,z_2)  K_h(x_2) \Big\}^2 \ d\Pb(z_1) \ d\Pb(z_2) \\
 \int \bigg[ \rho_h(x_1) &\{a_1 - \widehat\pi(x_1)\} \Big\{ {b}_h(x_1)^\T \widehat\Omega^{-1}   {b}_h(x_2) \Big\} \{y_2 - \widehat\mu_0(x_2)\} K_h(x_1) K_h(x_2) \bigg]^2 \ d\Pb(z_1) \ d\Pb(z_2) \\
%&\lesssim \frac{1}{h^{2d}} \int   \Big\{ {b}_h(x_1)^\T \widehat\Omega^{-1} {b}_h(x_2) \Big\}^2  K_h(x_1) K_h(x_2)   \ d\Pb(z_1) \ d\Pb(z_2) \\
&\lesssim \frac{1}{h^{2d}} \int \rho_h(x_1)^2 \Big\{ b_{hk}(x_1)^\T \widehat\Omega^{-1} b_{hk}(x_2) \Big\}^2  K_h(x_1) K_h(x_2) \ dF(x_1) \ dF(x_2) \\
&\lesssim \frac{1}{h^{2d}} \int  \Big\{ \overline{b}(v_1)^\T M \overline{b}(v_2) \Big\}^2 \ dF^*(v_1) \ dF^*(v_2) \\
%& =  \frac{1}{h^{2d}}\int \Big\{ \sum_{j,\ell} M_{j\ell} b^*_j(x_1) b^*_\ell(x_2) \Big\}^2 K_h(x_1) K_h(x_2) \ dF(x_1) \ dF(x_2) \\
%& =  \frac{1}{h^{2d}}\sum_{j,\ell}  \sum_{j',\ell'}  M_{j\ell} M_{j' \ell'} \int  b^*_j(x_1) b^*_{j'}(x_1)K_h(x_1) \ dF(x_1)  \int b^*_\ell(x_2)  b^*_{\ell'}(x_2)K_h(x_2)    \ dF(x_2) \\
%& = \frac{1}{h^{2d}} \sum_{j,\ell}  \sum_{j',\ell'}  M_{j\ell} M_{j' \ell'} \one(j=j')  \one(\ell=\ell')  \\
&=  \frac{\| M \|_2^2}{h^{2d}} 
%&= \|  \Omega^{1/2} \widehat\Omega^{-1} \Omega^{1/2} \|_2^2 / h^{2d} \\
= \left( \frac{1}{h^{2d}} \right)  \|  \Omega^{1/2} \Omega^{-1} \Omega^{1/2} + \Omega^{1/2} (\widehat\Omega^{-1} - \Omega^{-1}) \Omega^{1/2} \|_2^2 \\
& \leq 2\left( \frac{1}{h^{2d}} \right)  \left(  \| I \|_2^2 + \| \Omega^{1/2} (\widehat\Omega^{-1} - \Omega^{-1}) \Omega^{1/2} \|_2^2 \right) \\
%&= 2\left( \frac{1}{h^{2d}} \right)  \left( k + \| \Omega^{1/2} (\widehat\Omega^{-1} - \Omega^{-1}) \Omega^{1/2} \|_2^2 \right) \\
& \leq 2\left( \frac{1}{h^{2d}} \right) k \left( 1 + \| \Omega \|^2 \| \widehat\Omega^{-1} - \Omega^{-1} \|^2 \right)
\end{align*}
where the first line follows by definition,  the second since $(a-\widehat\pi)$ and $(y-\widehat\mu_0)$ are uniformly bounded and $0 \leq K_h(x) \lesssim h^{-d}$, the third by a change of variables $v=\frac{1}{2} + \frac{x-x_0}{h}$, by definition of $\overline{b}$ and $M$, and since $\rho(v)$ is bounded, the fourth by definition since $\int \overline{b} \overline{b}^\T \ dF^*=I$, and the last two by basic properties of the Frobenius norm (e.g., the triangle inequality and $\| A \|_2 \leq \sqrt{k} \| A \|$) together with  $(a + b)^2 \leq 2(a^2 + b^2)$. 
\end{proof}

\bigskip

\subsection{Technical Lemmas}

\begin{applemma} \label{lem:projections}
Let $\Pi_w f(t) = b(t)^\T \Omega^{-1} \int b(x) w(x) f(x) \ d\Pb(x)$ denote a $w$-weighted projection with $\Omega = \int b w b^T d\Pb$ for $w \geq 0$. And let $\theta_{w,f} = \Omega^{-1} \int b(x) w(x) f(x) \ d\Pb(x)$ denote the coefficients of the projection. Then:
\begin{enumerate}
\item[(i)] projections are orthogonal to residuals, i.e., 
$$ \int (\Pi_w f) (I-\Pi_w) g \ w d\Pb = 0 , $$
\item[(ii)] the $L_2(w\Pb)$ norm of a projection is no more than that of the function, i.e., 
$$ \int (\Pi_w f)^2 w \ d\Pb \leq \int f^2 \ w d\Pb , $$
and the $L_2(w\Pb)$ norm of the approximation error is no more than that of the function, i.e.,
$$ \int \{(I-\Pi_w) f\}^2 w \ d\Pb \leq \int f^2 \ w d\Pb , $$
\item[(iii)] the $L_2$ norm of the scaled coefficients is no more than the $L_2(w\Pb)$ norm of the function, i.e., 
$$ \| \Omega^{1/2} \theta_{w,f} \|^2 \leq \int f^2 \ wd\Pb . $$
\end{enumerate}
\end{applemma}

\medskip

\begin{proof}
For (i) let $b^*(x) = \Omega^{-1/2} b(x)$ so that $\int b^* w(b^*)^\T d\Pb = I$ and note that
\begin{align*}
\int ( \Pi_w f ) (I-\Pi_w) &g \ w d\Pb = \int \int b(x)^\T \Omega^{-1} b(t) f(t)g(x) \ wd\Pb(t) \ wd\Pb(x) \\
& \hspace{.9in} -  \int b(x)^\T \Omega^{-1} b(t) f(t) b(x)^\T \Omega^{-1}  b(u) g(u) \ wd\Pb(t) \ wd\Pb(x) \ wd\Pb(u) \\
&= \int b^*(x)^\T b^*(t) f(t) g(x) \ wd\Pb(t) \ wd\Pb(x) \\
& \hspace{.9in} -  \int f(t) b^*(t)^\T b^*(x) b^*(x)^\T  b^*(u) g(u) \ wd\Pb(t) \ wd\Pb(x) \ wd\Pb(u) \\
&= \int b^*(x)^\T b^*(t) f(t) g(x) \ wd\Pb(t) \ wd\Pb(x) \\
& \hspace{.9in} -  \int f(t) b^*(t)^\T  b^*(u) g(u) \ wd\Pb(t)  \ wd\Pb(u) = 0
\end{align*}
where the last line follows since $\int b^*(x) b^*(x)^\T \ wd\Pb(x) = I$. \\

For (ii) we have by definition that $\int \{\Pi_w f(x)\}^2 w(x) \ dF(x) $ equals
\begin{align}
 \int &\left\{ \int f(u) w(u) b(u)^\T \ d\Pb(u) \Omega^{-1} b(x) \right\} \left\{ b(x)^\T \Omega^{-1} \int b(t) w(t) f(t) \ d\Pb(t) \right\} w(x) d\Pb(x) \nonumber \\
&= \left\{ \int f(u) w(u) b(u)^\T \ d\Pb(u) \Omega^{-1} \right\} \left\{  \int b(t) w(t) f(t) \ d\Pb(t) \right\} \label{eq:l2normproj} \\
&= \int  f(t) \Pi_w f(t) w(t) \ d\Pb(t) \nonumber
\end{align}
and so
\begin{align*}
\int (f-\Pi_w f )^2 w \ d\Pb &= \int f^2 w \ d\Pb - 2 \int f \Pi_w f w \ d\Pb + \int (\Pi_w f)^2 w \ d\Pb  \\
&= \int f^2 w \ d\Pb -  \int (\Pi_w f)^2 w \ d\Pb 
\end{align*}
which implies the results, as long as $w \geq 0$, so that the far left and right sides are non-negative. \\

For (iii) we have
\begin{align*}
 \| \Omega^{1/2} \theta_{w,f} \|^2  &=  ( \Omega^{1/2} \theta_{w,f} )^\T (\Omega^{1/2} \theta_{w,f}) 
 =  \theta_{w,f}^\T \Omega \theta_{w,f} \\
 &= \left( \int b^\T w f \ d\Pb \right) \Omega^{-1}  \left( \int b w f \ d\Pb \right)  \\
 &= \int \{\Pi_w f(x)\}^2 w(x) \ d\Pb(x)
\end{align*}
where the last equality holds by \eqref{eq:l2normproj}, and so the result follows from Lemma \ref{lem:projections}(ii).
\end{proof}

\bigskip

\begin{applemma} \label{lem:eigen}
Assume:
\begin{enumerate}
\item $0 < b \leq \lambda_{\min}(\Omega) \leq \lambda_{\max}(\Omega) \leq B < \infty$
\item $0 < c \leq \| d \widehat{F}^*/dF^* \|_\infty \leq C < \infty$
\end{enumerate} 
Then
$$ bc \leq \lambda_{\min}(\widehat\Omega) \leq \lambda_{\max}(\widehat\Omega) \leq BC  $$
and
$$ \| \widehat\Omega^{-1} - \Omega^{-1} \| \leq  \left( \frac{ B }{b^2 c} \right)  \| (d\widehat{F}^*/dF^*) - 1 \|_\infty . $$
\end{applemma}

\bigskip

\begin{proof}
The logic mirrors that of the proof of Proposition 2.1 in  \citet{belloni2015some}. Note that
\begin{align*}
a^\T \widehat\Omega a &=  \int \Big\{ a^\T b(v) \Big\}^2 \ d\widehat{F}^*(v) \\
& \leq \| d\widehat{F}^*/dF^* \|_\infty  \int \Big\{ a^\T b(v) \Big\}^2  \ d{F}^*(v) \\
&= \| d\widehat{F}^*/dF^* \|_\infty \ a^\T \Omega a
\end{align*}
and by the same logic, $a^\T \Omega a \leq \| dF^*/d\widehat{F}^* \|_\infty \ a^\T \widehat\Omega a$. Therefore 
\begin{align*}
\lambda_{\max}( \widehat\Omega) &= \max_{\|a \|=1} a^\T \widehat\Omega a \leq  \| d\widehat{F}^*/dF^* \|_\infty \lambda_{\max}(\Omega) \\
\lambda_{\min}( \widehat\Omega) &=  \min_{\|a \|=1} a^\T \widehat\Omega a \geq \| dF^*/d\widehat{F}^* \|_\infty^{-1} \lambda_{\min}(\Omega)
\end{align*}
by the min-max theorem, which gives the first inequality.
For the second, note
\begin{align*}
 \| \widehat\Omega^{-1} - \Omega^{-1} \| &=  \| \widehat\Omega^{-1} (\Omega -\widehat\Omega ) \Omega^{-1} \| \\
 &\leq \| \widehat\Omega^{-1} \| \| \Omega -\widehat\Omega \| \| \Omega^{-1} \|
\end{align*}
by the sub-multiplicative property of the operator norm, and then
$$ \| \Omega -\widehat\Omega \|  \leq \| (d\widehat{F}^*/dF^*) - 1 \|_\infty \ \| \Omega \|  $$
by the same logic as above. 
\end{proof}

\bigskip

\subsection{Covariate Density}
\label{sec:appcovdensity}
 
Although the lower bounds in Theorems \ref{thm:lowerbound} and \ref{thm:lowerbound2} are given for a simple model where the covariate density is merely bounded, the rates are only attainable when the covariate density is known or can be estimated accurately enough, as discussed in Remark \ref{rem:cond1}. More specifically, the lower and upper bounds match for a model in which 
the covariate density 
$dF(x)$ is known, satisfies $\int \one\{ \|x-x_0 \| \leq h/2 \} \ dF(x) \asymp h^d$, and has local support $\{x \in \R^d: dF(x)>0, \|x-x_0\| \leq h/2  \}$ on a union of no more than $k$ disjoint cubes all with proportional volume, for $h$ and $k$ defined in Proposition \ref{prop:hellbound2}. The lower bounds come from Theorems \ref{thm:lowerbound} and \ref{thm:lowerbound2} since the density used in our construction satisfies these conditions, while the matching upper bound is implied by the results in Section \ref{sec:upperbound} together with the following lemma, in which we prove that Condition \ref{cnd:reg1} holds for this class of densities. \\

\begin{applemma} \label{lem:covdensity}
Assume:
\begin{enumerate}
\item $dF(x)$ satisfies $\int \one\{ \|x-x_0 \| \leq h/2 \} \ dF(x) \asymp h^d$, 
\item $dF(x)$ is bounded above and below away from zero on its local support $\{x \in \R^d: dF(x)>0, \|x-x_0\| \leq h/2  \}$, which is a union of no more than $k$ disjoint cubes all with proportional volume, for $h$ and $k$ defined in Proposition \ref{prop:hellbound2}.
\end{enumerate}
Let $dF^*(v)=dF(x_0 + h(v-1/2))$ denote the distribution in  $\mathcal{B}_h(x_0)$, the $h$-ball around $x_0$,  mapped to $[0,1]^d$, and similarly for $\pi^*$.  Then the eigenvalues of 
$$  Q %= \int \rho_h(x) K_h(x) \pi(x) \{1-\pi(x) \} \rho_h(x)^\T dF(x) $$
= \int  \rho(v)  \rho(v)^\T  \pi^*(v) \{ 1-\pi^*(v)\}\ dF^*(v) $$
are bounded above and below away from zero, and  there exists a basis $b$ with \Holder{} approximation property \eqref{eq:holder}, 
 for which the eigenvalues of
$$ \Omega = \int b(v) b(v)^\T \ dF^*(v) $$
are bounded above and below away from zero. 
\end{applemma}

\bigskip

\begin{proof}
Note that since $\epsilon(1-\epsilon) \leq \pi(1-\pi) \leq 1/4$ and since $dF(x)$ is bounded above and below on its support, the eigenvalues of $Q$ will be bounded if those of the matrix 
$$ \sum_{j=1}^k \int \rho(v) \rho(v)^\T \one(v \in M_j) \ dv $$
are, where $M_j$ indicates the $j$th disjoint cube making up the local support of $F$.  By the min-max theorem, the eigenvalues are bounded by the min and max of
$$  a^\T \int \rho(v) \rho(v)^\T \sum_j \one(v \in M_j) \ dv \ a = \int \Big\{ a^\T \rho(v) \Big\}^2 \sum_j \one(v \in M_j) \ dv   $$
over all $a \in \R^q$ with $\| a\|=1$. First consider lower bounding the eigenvalues. 
Note $g(v) = a^\T \rho(v)$ is a polynomial of degree at most $q$. Therefore
\begin{align*}
\sum_j \int g(v)^2 \one(v \in M_j) \ dv &\geq \int g(v)^2 \one\{ g(v)^2 \geq \epsilon \} \sum_j \one(v \in M_j) \ dv \\
& \geq \epsilon \int \one\{ g(v)^2 \geq \epsilon \} \sum_j \one(v \in M_j) \ dv \\
& \geq \epsilon \left\{ \sum_j \int \one(v \in M_j) \ dv - \int \one\{ g(v)^2 < \epsilon \} \ dv \right\} \\
& \geq \epsilon \left( C^* - C \epsilon^{1/2q} \right)  = \left( \frac{C^*}{\sqrt{C}} \right)^{4q} > 0 
\end{align*}
where the last line follows since
\begin{align*}
\int \sum_j \one( v \in M_j) \ dv & \asymp \int dF^*(v) = \int dF(x_0 + h(v-1/2)) \\
&= h^{-d} \int \one\{ \|x-x_0 \| \leq h/2 \} \ dF(x) \geq C^*
\end{align*}
from a change of variables $x = x_0 + h(v-1/2)$, so that $v = 1/2 + (x-x_0)/h$ and $h^d \ dv = dx$, together with Assumption 1 of the lemma,  and by Theorem 4 of \citet{carbery2001distributional}. The last equality follows if we choose $\epsilon =\left(  C^*/2C \right)^{2q}$. To upper-bound the eigenvalues, note that
\begin{align*}
 \int \Big\{ a^\T \rho(v) \Big\}^2 \sum_j \one(v \in M_j) \ dv  \leq  \int \Big\{ a^\T \rho(v) \Big\}^2  \ dv = \| a \|^2 = 1
\end{align*}
since $\int \rho(v) \rho(v)^\T \ dv = I$ by the orthonormality of $\rho$. \\

Now consider the eigenvalues of $\Omega$. We will construct a basis of order $k$ for which eigenvalues of $\Omega$ are bounded and for which the \Holder{} approximation property \eqref{eq:holder} holds.  \\

\begin{remark}
For simplicity we only consider the case where there are exactly $k$ cubes in the local support of $F$; if there are finitely many cubes, say $M < \infty$, the arguments are more straightforward, orthonormalizing a standard series $M$ times, once per cube, and taking $k/M$ terms from each. We omit details in the intermediate case where the number of cubes scales at a rate slower than $k$, but we expect similar arguments to those below can be used. \\
\end{remark}

First, let $\rho(x)$ denote the tensor products of a Legendre basis of order $s$, orthonormal on $[-1,1]$, i.e.,  tensor products of
$$ \rho_j(x_\ell) = \frac{\sqrt{ j+1/2} }{2^j j!}  \frac{d^j}{dx_\ell^j} (x_\ell^2-1)^j  $$
for $j \in \{1,\dots,s\}$ and $\ell \in \{1,\dots,d\}$. This basis satisfies
$$ \int \rho(x) \rho(x)^\T \ dx = I_s . $$
Now shift and rescale so that the basis is orthonormal on $M_j$ with
% for [0,1]: subtract off left edge & multiply by 1/length
% for [-1,1]: subtract off midpoint & multiply by 2/length
$$ b_j(x) = 2^d \sqrt{ k} \rho\left(4k^{1/d} \left(x- \frac{1}{2} - \frac{m_j - x_0}{h} \right) \right) $$
and so satisfies
\begin{align*}
\int b_j(x) b_j(x)^\T \one(x \in M_j) \ dx &=4^d k \int \rho\left(4k^{1/d} \left(x- \frac{1}{2} - \frac{m_j - x_0}{h} \right) \right) \one(x \in M_j) \ dx \\
&= \int_{-1}^1 \rho(v) \rho(v)^\T \ dv = I_s
\end{align*}
where we used the change of variable $v = 4k^{1/d} \left(x- \frac{1}{2} - \frac{m_j - x_0}{h} \right)$ so that $dv = 4^d k \ dx$ and $x = v/4k^{1/d} + 1/2 + (m_j-x_0)/h \in M_j$ when $v \in [-1,1]^d$.  Finally define the basis 
$$ b(v) = \{ b_1(v)^\T \one(v \in M_1), \dots, b_k(v)^\T \one(v \in M_k) \}^\T  $$
of length $sk$. Then
\begin{align*}
b(v) b(v)^\T = 
\begin{pmatrix} 
b_1(v) b_1(v)^\T \one(v \in M_1) & 0 & \cdots & 0 \\
0 & b_2(v) b_2(v)^\T \one(v \in M_2) & \cdots & 0 \\
\vdots & \vdots & \ddots & \vdots \\
0 & 0 & \cdots & b_k(v) b_k(v)^\T \one(v \in M_k) \\
\end{pmatrix}
\end{align*}
Therefore since $dF^*(x)$ is bounded above and below on its support, the $j$-th diagonal block of $\Omega$ is proportional to
\begin{align*}
\int   b_j(v) b_j(v)^\T \one(v \in M_j) \ dv&= I_s .
\end{align*}
Therefore the eigenvalues of $\Omega$ are all proportional to one and bounded as desired. Further, by  the same higher-order kernel arguments for local polynomials as in Proposition \ref{prop:projbias}, the basis satisfies the  \Holder{} approximation property \eqref{eq:holder}. 
\end{proof}

%%%%%%%%%%%%
%%%%%%%%%%%%
%%%%%%%%%%%%
%%%%%%%%%%%%
%%%%%%%%%%%%

\end{document}